\input ./style/arxiv-general.cfg
\documentclass[aos,MSNbibl,seceqn,nameyear,dvips]{arximspdf}
\makeatletter
   \@ifpackageloaded{graphicx}{}{\usepackage{graphicx}}
\makeatother
\usepackage{url,breakurl}

\doi{10.1214/14-AOS1301}
\volume{43}
\issue{3}
\pubyear{2015}
\firstpage{1117}
\lastpage{1140}
\docsubty{FLA}

\makeatletter
\def\sfrac#1#2{#1/#2}
\def\vfrac#1#2{(#1)/#2}

\newcommand{\Cov}{\operatorname{Cov}}
\newcommand{\rrvert}{\vert}
\newcommand{\rrVert}{\Vert}
\newcommand{\llvert}{\vert}
\newcommand{\llVert}{\Vert}
\newtheorem{theorem}{Theorem}[section]
\newtheorem{corollary}{Corollary}[section]
\newproclaim{remark}{Remark}[section]
\newcommand{\X}{\mathbf{X}}
\renewcommand{\S}{\mathbf{S}}
\newcommand{\D}{\mathbf{D}}
\newcommand{\Y}{\mathbf{Y}}
\newcommand{\V}{\mathbf{V}}
\renewcommand{\L}{\mathbf{L}}
\newcommand{\N}{\mathbb{N}}
\newcommand{\R}{\mathbb{R}}
\newcommand{\Z}{\mathbb{Z}}
\newcommand{\C}{\mathbf{C}}
\newcommand{\overset}{\stackrel}
\makeatother

\begin{document}
\begin{frontmatter}

\title{Covariance matrix estimation and linear process bootstrap for
multivariate time series of possibly increasing dimension}
\runtitle{Covariance matrix estimation and multivariate LPB}

\begin{aug}
\author[A]{\fnms{Carsten}~\snm{Jentsch}\ead[label=e1]{cjentsch@mail.uni-mannheim.de}\thanksref{T1}}
\and
\author[B]{\fnms{Dimitris N.}~\snm{Politis}\corref{}\ead[label=e2]{politis@math.ucsd.edu}\thanksref{T2}}
\runauthor{C. Jentsch and D. N. Politis}
\affiliation{University of Mannheim and University of California, San Diego}
\address[A]{Department of Economics\\
University of Mannheim\\
L7, 3-5\\
68131 Mannheim\\
Germany\\
\printead{e1}}
\address[B]{Department of Mathematics\\
University of California, San Diego\\
La Jolla, California 92093-0112\\
USA\\
\printead{e2}}
\end{aug}
\thankstext{T1}{Supported by a fellowship within the Ph.D. Program of the German Academic Exchange Service (DAAD) when
Carsten Jentsch was visiting the University of California, San Diego, USA.}
\thankstext{T2}{Supported in part by NSF
Grants DMS-13-08319 and DMS-12-23137.}

%
\received{\smonth{1} \syear{2014}}
%
\revised{\smonth{9} \syear{2014}}

%
\begin{abstract}
Multivariate time series present many challenges, especially when they
are high dimensional. The paper's focus is twofold. First, we
address the subject of consistently estimating the autocovariance
sequence; this is a sequence of matrices that we conveniently
stack into one huge matrix. We are then able to show consistency of an
estimator based on the so-called
\textit{flat-top tapers}; most importantly, the consistency holds true
even when the time series dimension is allowed to increase with the
sample size.
Second, we revisit the linear process bootstrap (LPB) procedure
proposed by
McMurry and Politis [\textit{J.~Time Series Anal.} \textbf{31} (2010) 471--482]
for univariate time series. Based on the aforementioned stacked
autocovariance matrix estimator, we are able to define a version of the
LPB that is valid for
multivariate time series. Under rather general assumptions, we show
that our multivariate linear process bootstrap (MLPB) has asymptotic validity
for the sample mean in two important cases: (a) when the time series
dimension is fixed and (b) when it is allowed to increase with sample size.
As an aside, in case (a) we show that the MLPB works also for spectral
density estimators which is a novel result even in the univariate case.
We conclude with a simulation study that demonstrates the superiority
of the MLPB in some important cases.
\end{abstract}

%
\end{frontmatter}

\section{Introduction}\label{introduction}
Resampling methods for dependent data such as time series have been studied
extensively over the last decades. For an overview of existing
bootstrap methods see the monograph of \citet{Lahiri2003} and the review
papers by \citet{Buehlmann}, \citet{Paparoditis2002}, \citet{HaerdleHorowitzKreiss2003}, \citet{Politis2003a} or the recent review paper by \citet{KreissPaparoditis2011}.
Among the most popular bootstrap procedures in time series analysis, we
mention the autoregressive (AR) sieve bootstrap [cf. \citeauthor{Kreiss1992} (\citeyear{Kreiss1992,Kreiss1999}),
\citet{Buehlmann1997}, \citet{KreissPaparoditisPolitis2011}] and
block bootstrap and its variations; cf. \citet{Kuensch1989}, \citet{LiuSingh1992},
\citeauthor{PolitisRomano1992} (\citeyear{PolitisRomano1992,PolitisRomano1994}), etc.
A recent addition to the available time series bootstrap methods was
the linear process bootstrap (LPB) introduced by \citet{McMurryPolitis2010} who showed its validity for the sample mean for univariate
stationary processes without actually assuming linearity of the
underlying process.

The main idea of the LPB is to consider the time series data of length
$n$ as one large $n$-dimensional vector and to estimate appropriately
the entire covariance structure of this vector. This is executed by
using tapered covariance matrix estimators based on flat-top kernels
that were defined in \citet{Politis2001}. The resulting covariance matrix is
used to whiten the data by pre-multiplying the original (centered) data
with its inverse Cholesky matrix; a modification of the eigenvalues, if
necessary, ensures positive definiteness. This decorrelation property
is illustrated in Figures~5 and 6 in \citet{JentschPolitis2013}. After
suitable centering and standardizing, the whitened vector is treated as
having independent and identically
distributed (i.i.d.) components with zero mean and unit variance.
Finally, i.i.d. resampling from this vector and pre-multiplying the
corresponding bootstrap vector of residuals with the Cholesky matrix
itself results in a bootstrap sample that has (approximately) the same
covariance structure as the original time series.

Due to the use of flat-top kernels with compact support, an abruptly
dying-out autocovariance structure is induced to the bootstrap
residuals. Therefore, the LPB is particularly suitable for---but not
limited to---time series of moving average (MA) type. In a sense, the
LPB could be considered the closest
analog to an MA-sieve bootstrap which is not practically feasible due
to nonlinearities in the estimation of the MA parameters.
A further similarity of the LPB to MA fitting, at least in the
univariate case, is the
equivalence of computing the Cholesky decomposition of the covariance
matrix to the innovations algorithm; cf. \citet{RissanenBarbosa1969},
\citet{BrockwellDavis1988} and \citet{BrockwellMitchell1997}, the
latter addressing the multivariate case.

Typically, bootstrap methods extend easily from the univariate to the
multivariate case, and the same is true for time series
bootstrap procedures such as the aforementioned AR-sieve bootstrap
and the block bootstrap. By contrast, it has not been clear to date
if/how the LPB could be
successfully applied in the context of multivariate time series data;
a proposal to that effect was described in \citet{JentschPolitis2013}---who refer to an earlier
preprint of the paper at hand---but
it has been unclear to date whether the multivariate LPB is
asymptotically consistent
and/or if it competes well with other methods.
Here we attempt to fill this gap: we show how to implement the LPB in a
multivariate context and prove its validity for the sample mean and for
spectral density estimators, the latter being a new result even in the
univariate case.
{Note that the limiting distributions of the sample mean and of kernel
spectral density estimators depend only on the second-order moment
structure. Hence it is intuitive that the LPB would be well suited for
such statistics since it generates a linear process in the bootstrap
world that mimics well the second-order moment structure of the real
world.} Furthermore, in the spirit of the times, we consider the
possibility that the time series dimension is increasing with sample
size and
identify conditions under which the multivariate linear process
bootstrap (MLPB)
maintains its asymptotic validity, even in this case.
The key here is to address the subject of consistently estimating the
autocovariance sequence; this is a sequence of matrices that we conveniently
stack into one huge matrix. We are then able to show consistency of an
estimator based on the aforementioned flat-top tapers; most
importantly, the consistency holds true even when the time series
dimension is allowed to increase with the sample size.

The paper is organized as follows. In Section~\ref{preliminaries}, we
introduce the notation of this paper, discuss tapered covariance matrix
estimation for multivariate stationary time series and state
assumptions used throughout the paper; we then present our results on
convergence with respect to operator norm of tapered covariance matrix
estimators. The MLPB bootstrap algorithm and some remarks can be found
in Section~\ref{bootstrapscheme}, and results concerned with validity
of the MLPB for the sample mean and kernel spectral density estimates
are summarized in Section~\ref{asymptoticresults}. Asymptotic results
established for the case of increasing time series dimension are stated
in Section~\ref{asymptoticresultsincreasing}, where operator norm
consistency of tapered covariance matrix estimates and a validity
result for the sample mean are discussed. A finite-sample simulation
study is presented in Section~\ref{secsim}.
Finally, all proofs, some additional simulations and a real data
example on the weighted mean of an increasing number of stock prices
taken from the German stock index DAX can be found at the paper's
supplementary material [\citet{JentschPolitis2014}], which is also available
at \url{http://www.math.ucsd.edu/\textasciitilde politis/PAPER/MLPBsupplement.pdf}.

\section{Preliminaries}\label{preliminaries}

Suppose we consider an $\R^d$-valued time series process $\{\underline
{X}_t, t\in\Z\}$ with $\underline{X}_t=(X_{1,t},\ldots,X_{d,t})^T$,
and we have data $\underline{X}_1,\ldots,\underline{X}_n$ at hand.
The process $\{\underline{X}_t, t\in\Z\}$ is assumed to be strictly
stationary and its $(d\times d)$ autocovariance matrix $\C
(h)=(C_{ij}(h))_{i,j=1,\ldots,d}$ at lag $h\in\Z$ is
\begin{equation}
\C(h)=E \bigl((\underline{X}_{t+h}-\underline{\mu}) (\underline
{X}_t-\underline{\mu})^T \bigr), \label{covariance}
\end{equation}
where $\underline{\mu}=E(\underline{X}_t)$, and the sample
autocovariance $\widehat\C(h)=(\widehat C_{ij}(h))_{i,j=1,\ldots,d}$
at lag $\llvert  h\rrvert <n$ is defined by
\begin{equation}
\widehat\C(h)=\frac{1}{n}\sum_{t=\max(1,1-h)}^{\min(n,n-h)}
(\underline{X}_{t+h}-\overline{\underline{X}}) (\underline
{X}_t-\overline{\underline{X}})^T, \label{samplecovariance}
\end{equation}
where $\overline{\underline{X}}=\frac{1}{n}\sum_{t=1}^n \underline
{X}_t$ is the $d$-variate sample mean vector. Here and throughout the
paper, all matrix-valued quantities are written as bold letters, all
vector-valued quantities are underlined, $\mathbf{A}^T$ indicates the
transpose of a matrix $\mathbf{A}$, $\overline{\mathbf{A}}$ the
complex conjugate of $\mathbf{A}$ and $\mathbf{A}^H=\overline
{\mathbf{A}}^T$ denotes the transposed conjugate of $\mathbf{A}$.
Note that it is also possible to use unbiased sample autocovariances,
that is, having $n-\llvert  h\rrvert $ instead of $n$ in the
denominator of (\ref
{samplecovariance}). Usually the biased version as defined in (\ref
{samplecovariance}) is preferred because it guarantees a positive
semi-definite estimated autocovariance function, but our tapered
covariance matrix estimator discussed in Section~\ref{taperedestimator} is
adjusted in order to become positive definite in any case.

Now, let $\underline{X}=\operatorname{vec}(\X)=(X_1,\ldots,X_{dn})^T$ be the
$dn$-dimensional vectorized version of the $(d\times n)$ data matrix
$\X=[\underline{X}_1\dvtx \underline{X}_2\dvtx \cdots\dvtx \underline{X}_n]$, and
denote the covariance matrix of $\underline{X}$, which is symmetric
block Toeplitz, by $\bolds{\Gamma}_{dn}$, that is,
\begin{equation}
\bolds{\Gamma}_{dn}= \pmatrix{\C(i-j)
\cr
i,j=1,\ldots,n}= \pmatrix{
\Gamma_{dn}(i,j)
\cr
i,j=1,\ldots,dn}, \label{Gammadn}
\end{equation}
where $\Gamma_{dn}(i,j)=\Cov(X_i,X_j)$ is the covariance between the
$i$th and $j$th entry of~$\underline{X}$. Note that the second order
stationarity of $\{\underline{X}_t,t\in\Z\}$ does \emph{not} imply
second-order stationary behavior of the vectorized $dn$-dimensional
data sequence $\underline{X}$. This means that the covariances $\Gamma
_{dn}(i,j)$ truly depend on both $i$ and $j$ and not only on the
difference $i-j$. However, the following one-to-one correspondence
between $\{C_{ij}(h),h\in\Z,i,j=1,\ldots,d\}$ and $\{\Gamma
_{dn}(i,j),i,j\in\Z\}$ holds true. Precisely, we have
\begin{eqnarray}
\Gamma_{dn}(i,j) &=& \Cov(X_i,X_j)
\nonumber
\\
&=& \Cov(X_{m_1(i),m_2(i)},X_{m_1(j),m_2(j)}) \label{relation}
\\
&=& C_{\underline{m}_1(i,j)}\bigl(m_2(i,j)\bigr),
\nonumber
\end{eqnarray}
where $\underline{m}_1(i,j)=(m_1(i),m_1(j))$ and
$m_2(i,j)=m_2(i)-m_2(j)$ with $m_1(k)=(k-1)\operatorname{mod} d+1$ and $m_2(k)=
\lceil k/d \rceil$, and $\lceil x\rceil$ denotes the smallest
integer greater or equal to $x\in\R$.

If one is interested in estimating the quantity $\bolds{\Gamma}_{dn}$, it seems natural
to plug in the sample covariances $\widehat{\mathbf C}(i-j)$ and $\widehat\Gamma_{dn}(i,j)=\widehat
C_{\underline{m}_1(i,j)}(m_2(i,j))$ in $\bolds{\Gamma}_{dn}$ and to use
\[
\widehat{\bolds \Gamma}_{dn}= \pmatrix{ \widehat\C(i-j)
\cr
i,j=1,\ldots,n }= \pmatrix{ \widehat\Gamma_{dn}(i,j)
\cr
i,j=1,\ldots,dn}.
\]
But unfortunately this estimator is \emph{not} a consistent estimator
for $\bolds{\Gamma}_{dn}$ in the sense that the operator norm of
$\widehat{\bolds \Gamma}_{dn}-\bolds{\Gamma}_{dn}$ does not
converge to zero.
This was shown by \citet{WuPourahmadi2009}, and to dissolve this
problem in the univariate case, they proposed a banded estimator of the
sample covariance matrix to achieve consistency.
This has been generalized by \citet{McMurryPolitis2010}, who considered
general flat-top kernels as weight functions.

In Section~\ref{taperedestimator}, we follow the paper of \citet{McMurryPolitis2010} and propose a tapered estimator of $\bolds{\Gamma
}_{dn}$ and show its consistency in Theorem~\ref
{operatornormconvergence1} for the case of multivariate processes.
Moreover, we state a modified estimator that is guaranteed to be
positive definite for any finite sample size and show its consistency
in Theorem~\ref{operatornormconvergence2} and of related quantities in
Corollary~\ref{operatornormconvergence3}. But prior to this, we state
the assumptions that are used throughout this paper in the following.

\subsection{Assumptions}\label{assumptions}

\begin{longlist}[(A2)]
\item[(A1)] $\{\underline{X}_t,t\in\Z\}$ is an $\R^d$-valued
strictly stationary time series process with mean $E(\underline
{X}_t)=\underline{\mu}$ and
autocovariances $\C(h)$ defined in (\ref{covariance}) such that $\sum_{h=-\infty}^\infty\llvert  h\rrvert ^g \llvert \mathbf
{C}(h)\rrvert _1<\infty$ for some $g\geq
0$ to be further specified.
Let $\llvert \mathbf{A}\rrvert _p=(\sum_{i,j} \llvert
a_{ij}\rrvert ^p)^{1/p}$ for some matrix
$\mathbf{A}=(a_{ij})$.
\item[(A2)] There exists a constant $M<\infty$ such that for all
$n\in\N$, all $h$ with $\llvert  h\rrvert <n$ and all
$i,j=1,\ldots,d$, we have
\[
\Biggl\llVert \sum_{t=1}^{n}(X_{i,t+h}-
\overline{X}_i) (X_{j,t}-\overline {X}_j)-nC_{ij}(h)
\Biggr\rrVert _2\leq M\sqrt{n},
\]
where $\llVert \mathbf{A}\rrVert _p=(E(\llvert \mathbf{A}\rrvert _p^p))^{1/p}$.
\item[(A3)] There exists an $n_0\in\N$ large enough such that for
all $n\geq n_0$ the eigenvalues $\lambda_1,\ldots,\lambda_{dn}$ of
the $(dn\times dn)$ covariance matrix $\bolds{\Gamma}_{dn}$ are
bounded uniformly away from zero.
\item[(A4)] {Define the projection operator $P_k(\underline
{X})=E(\underline{X}\llvert \mathcal{F}_{k})-E(\underline{X}\rrvert \mathcal
{F}_{k-1})$ for $\mathcal{F}_k=\sigma(\underline{X}_t,t\leq k)$, and
suppose that for all $i=1,\ldots,d$, we have $\sum_{m=0}^\infty\llVert
P_0\times\break X_{i,m}\rrVert _q<\infty$ and $\llVert \overline{X}_i-\mu
_i\rrVert _q=O(\frac
{1}{\sqrt{n}})$, respectively, for some $q\geq2$ to be further specified.}
\item[(A5)] For the sample mean, a CLT holds true. That is, we have
\[
\sqrt{n}(\overline{\underline{X}}-\underline{\mu})\overset {\mathcal{D}} {
\longrightarrow}\mathcal{N}(\underline{0},\mathbf{V}),
\]
where $\overset{\mathcal{D}}{\longrightarrow}$ denotes weak
convergence, $\mathcal{N}(\underline{0},\mathbf{V})$ is a normal
distribution with zero mean vector and covariance matrix $\mathbf
{V}=\sum_{h=-\infty}^\infty\C(h)$ with $\mathbf{V}$ positive definite.\vspace*{2pt}
\item[(A6)] For kernel spectral density estimates $\widehat
f_{pq}(\omega)$ as defined in (\ref{fhat}) in Section~\ref{asymptoticresults}, a CLT holds true. {That is, for arbitrary
frequencies $0\leq\omega_1,\ldots,\omega_s\leq\pi$}, we have that
\[
\sqrt{nb} \bigl(\widehat f_{pq}(\omega_j)-f_{pq}(
\omega _j)\dvtx p,q=1,\ldots,d; j=1,\ldots,s \bigr)
\]
converges\vspace*{1pt} to an $sd^2$-dimensional normal distribution for
$b\rightarrow0$ and $nb\rightarrow\infty$ such that $nb^5=O(1)$ as
$n\rightarrow\infty$, where the limiting covariance matrix is
obtained from
\begin{eqnarray*}
&& nb\Cov \bigl(\widehat f_{pq}(\omega),\widehat f_{rs}(\lambda)
\bigr)
\\
&&\qquad = \bigl(f_{pr}(\omega)\overline{f_{qs}(\omega)}
\delta_{\omega
,\lambda}+f_{ps}(\omega)\overline{f_{qr}(\omega)}
\tau_{0,\pi
} \bigr)
\frac{1}{2\pi}\int K^2(u)\,du+o(1)
\end{eqnarray*}
and the limiting bias from
\[
E \bigl(\widehat f_{pq}(\omega) \bigr)-f_{pq}(
\omega)=b^2 f''_{pq}(\omega)
\frac{1}{4\pi}\int K(u)u^2\,du+o \bigl(b^2 \bigr)
\]
for all $p,q,r,s=1,\ldots,d$, where $\delta_{\omega,\lambda}=1$ if
$\omega=\lambda$ and $\tau_{0,\pi}=1$ if $\omega=\lambda\in\{
0,\pi\}$ and zero otherwise, respectively. Therefore, $\mathbf
{f}(\omega)$ is assumed to be component-wise twice differentiable with
Lipschitz-continuous second derivatives.
\end{longlist}

Assumption \textup{(A1)} is quite standard, and the uniform convergence of
sample autocovariances in (A2) is satisfied under different types of
conditions (cf. \mbox{Remark}~\ref{remarkA2} below) and appears to be a
crucial condition here.
The uniform boundedness of all eigenvalues away from zero in (A3) is
implied by a nonsingular spectral density matrix $\mathbf{f}$ of
$(\underline{X}_t,t\in\Z)$. This follows with (\ref{Gammadn}) and
the inversion formula~from
\[
{\underline{c}^T\bolds{\Gamma}_{dn}\underline {c}=
\underline{c}^T \biggl(\int_{-\pi}^\pi
\mathbf{J}_\omega ^T\mathbf{f}(\omega)\overline{
\mathbf{J}}_\omega \,d\omega \biggr)\underline{c}\geq2\pi\llvert
\underline{c}\rrvert _2^2 \inf_\omega
\lambda _{\min}\bigl(\mathbf{f}(\omega)\bigr)}
\]
for all $\underline{c}\in\R^{dn}$, where $\mathbf{J}_\omega
=(e^{-i1\omega},\ldots,e^{-in\omega})\otimes\mathbf{I}_d$ and
$\otimes$ denotes the Kronecker product. 
The requirement of condition (A3) fits into the theory for the
univariate autoregressive sieve bootstrap as obtained in \citet{KreissPaparoditisPolitis2011}. Similarly,
a nonsingular spectral density matrix $\mathbf{f}$ implies positive
definiteness of the long-run variance $\mathbf{V}=2\pi\mathbf{f}(0)$
defined in (A5). Assumption \textup{(A4)} is, for instance, fulfilled if the
underlying process is linear or $\alpha$-mixing with summable mixing
coefficients by Ibragimov's inequality; cf., for example, \citet{Davidson1994}, Theorem 14.2. To achieve validity of the MLPB for the sample
mean and for kernel spectral density estimates in Section~\ref{asymptoticresults}, we have to assume unconditional CLTs in (A5) and
(A6), which are satisfied also under certain mixing conditions [cf.
\citet{Doukhan1994}, \citet{Brillinger1981}], linearity [cf. \citet{BrockwellDavis1991}, \citet{Hannan1970}] or weak dependence [cf. \citet{Dedeckeretal2007}].
Note also that the condition $nb^5=O(1)$ includes the optimal bandwidth
choice $nb^5\rightarrow C^2$, $C>0$ for second-order kernels, which
leads to a nonvanishing bias in the limiting normal distribution.

%
\begin{remark}\label{remarkA2}
Assumption (A2) is implied by different types of conditions imposed on
the underlying process $\{\underline{X}_t,t\in\Z\}$. We present
sufficient conditions for (A2) under linearity, mixing or weak dependence type conditions.
More precisely, (A2) is satisfied
if the process $\{\underline{X}_t,t\in\Z\}$ fulfills one of the
following conditions:
\begin{longlist}[(iii)]
\item[(i)] \textit{Linearity}. Suppose the process is linear; that
is, $\underline{X}_t=\sum_{k=-\infty}^\infty\mathbf{B}_k\underline
{e}_{t-k}$, $t\in\Z$, where $\{\underline{e}_t,t\in\Z\}$ is an
i.i.d. white noise with finite fourth moments
$E(e_{i,t}e_{j,t}e_{k,t}e_{l,t})<\infty$ for all $i,j,k,l=1,\ldots,d$
and the sequence of $(d\times d)$ coefficient matrices $\{\mathbf
{B}_k,k\in\Z\}$ is component-wise absolutely summable.
\item[(ii)] \textit{Mixing-type condition}. Let
\[
\operatorname{cum}_{a_1,\ldots,a_k}(u_1,\ldots,u_{k-1})=\operatorname{cum}(X_{a_1,u_1},
\ldots ,X_{a_{k-1},u_{k-1}},X_{a_k,0})
\]
denote the $k$th order joint cumulant of $X_{a_1,u_1},\ldots
,X_{a_{k-1},u_{k-1}},X_{a_k,0}$ [cf.\break \citet{Brillinger1981}], and suppose
$\sum_{s,h=-\infty}^\infty\llvert  \operatorname{cum}_{i,j,i,j}(s+h,s,h)\rrvert
<\infty$ for all
$i,j=1,\ldots,d$. Note that this is satisfied if $\{\underline
{X}_t,t\in\Z\}$ is $\alpha$-mixing such that $E(\llvert  X_1\rrvert _2^{4+\delta
})<\infty$ and $\sum_{j=1}^\infty j^2\alpha(j)^{\delta/(4+\delta
)}<\infty$ for some $\delta>0$; cf. \citet{Shao2010}, page 221.
%
%
\item[(iii)] \textit{Weak dependence-type condition}. Suppose for
all $i,j=1,\ldots,d$, we have
\[
\bigl\llvert \Cov\bigl((X_{i,t+h}-\mu_i) (X_{j,t}-
\mu _j),(X_{i,t+s+h}-\mu_i) (X_{j,t+s}-
\mu_j)\bigr)\bigr\rrvert \leq C\cdot\nu_{s,h},
\]
where $C<\infty$ and $(\nu_{s,h},s,h\in\Z)$ is absolutely summable,
that is,
\[
\sum_{s,h=-\infty}^\infty\llvert \nu_{s,h}\rrvert
<\infty;
\]
cf.
\citet{Dedeckeretal2007}.
\end{longlist}
\end{remark}

\subsection{Tapered covariance matrix estimation of multiple time series data}\label{taperedestimator}
To adopt the technique of \citet{McMurryPolitis2010}, let
\begin{eqnarray}
\kappa(x)=\cases{ 1, &\quad$\llvert x\rrvert \leq1$,
\cr
0, &\quad$\llvert x
\rrvert >c_\kappa$,
\cr
g\bigl(\llvert x\rrvert \bigr), &\quad
otherwise}\label{kappa}
\end{eqnarray}
be a so-called \textit{flat-top taper} [cf. \citet{Politis2001}], where
$\llvert  g(x)\rrvert <1$ and $c_\kappa\geq1$. The $l$-scaled
version of $\kappa
(\cdot)$ is defined by $\kappa_l(x)=\kappa (\frac{x}{l} )$
for some $l>0$.
As \citet{Politis2011} argues, it is advantageous to have a smooth taper
$\kappa$, so
the truncated kernel that corresponds to $g(x)=0$ for all $x$ is not
recommended.
The simplest example of a continuous taper function $\kappa$ with
$c_\kappa>1$ is the trapezoid
\begin{eqnarray}
\label{trapezoid} \kappa(x)=\cases{ 1, &\quad$\llvert x\rrvert \leq1$,
\cr
2-\llvert
x\rrvert , &\quad$1<\llvert x\rrvert \leq2$,
\cr
0, &\quad$\llvert x\rrvert >2$}
\end{eqnarray}
which is used in Section~\ref{secsim} for the simulation study; the
trapezoidal taper was first proposed by \citet{PolitisRomano1995} in
a spectral estimation setup. Observe also that the banding parameter
does not need to be an integer. The tapered estimator $\widehat{\bolds\Gamma}_{\kappa,l}$ of $\bolds{\Gamma}_{dn}$ is given by
\begin{equation}
\widehat{\bolds \Gamma}_{\kappa,l}= \pmatrix{\kappa_l(i-j)
\widehat\C(i-j)
\cr
i,j=1,\ldots,n}= \pmatrix{\widehat\Gamma_{\kappa,l}(i,j)
\cr
i,j=1,\ldots,dn}, \label{gammahat}
\end{equation}
where $\widehat\Gamma_{\kappa,l}(i,j)=\widehat C_{\underline
{m}_1(i,j)}^{\kappa,l} (m_2(i,j) )$ and $\widehat
C_{i,j}^{\kappa,l} (h )=\kappa_{l}(h)\widehat C_{i,j}(h)$.

The following Theorem~\ref{operatornormconvergence1} deals with
consistency of the tapered estimator $\widehat{\bolds \Gamma
}_{\kappa,l}$ with respect to operator norm convergence. It extends
Theorem 1 in \citet{McMurryPolitis2010} to the multivariate case and
does not rely on the concept of physical dependence only. The operator
norm of a complex-valued $(d\times d)$ matrix $\mathbf{A}$ is defined by
\[
\rho(\mathbf{A})=\max_{\underline{x}\in\mathbb{C}^d\dvtx \llvert
\underline
{x}\rrvert _2=1}\llvert \mathbf{A}\underline{x}
\rrvert _2,
\]
%
and it is well known that $\rho^2(\mathbf{A})=\lambda_{\max
}(\mathbf{A}^H\mathbf{A})=\lambda_{\max}(\mathbf{A}\mathbf
{A}^H)$, where $\lambda_{\max}(\mathbf{B})$ denotes the largest
eigenvalue of a matrix $\mathbf{B}$; cf. \citet{HornJohnson1990}, page 296.

%
\begin{theorem}\label{operatornormconvergence1}
Suppose that assumptions \textup{(A1)} with $g=0$ and \textup{(A2)} are satisfied. Then
it holds
\begin{eqnarray}\label{convergenceeq}
&& \bigl\llVert \rho(\widehat{\bolds \Gamma}_{\kappa,l}-\bolds{\Gamma
}_{dn})\bigr\rrVert _2
\nonumber\\[-8pt]\\[-8pt]\nonumber
&&\qquad \leq \frac{4Md^2(\lfloor c_\kappa l\rfloor+1)}{\sqrt{n}}+2\sum
_{h=0}^{\lfloor c_\kappa l\rfloor}\frac{\llvert  h\rrvert
}{n}\bigl\llvert
\mathbf{C}(h)\bigr\rrvert _1
+2\sum_{h=l+1}^{n-1}\bigl\llvert
\mathbf{C}(h)\bigr\rrvert _1.
\end{eqnarray}
%
\end{theorem}

The second term on the right-hand side of (\ref{convergenceeq}) can
be represented as $2\frac{\lfloor c_\kappa l\rfloor}{n}\sum_{h=0}^{\lfloor c_\kappa l\rfloor}\frac{\llvert  h\rrvert
}{\lfloor c_\kappa
l\rfloor}\llvert \mathbf{C}(h)\rrvert _1$ and vanishes
asymptotically due to the
Kronecker\vadjust{\goodbreak} lemma and is of order $o(l/n)$. The
third one converges to zero for $l=l(n)\rightarrow\infty$ as
$n\rightarrow\infty$ and the leading first term for $l=o(\sqrt{n})$.
Hence, for
$1/l+l/\sqrt{n}=o(1)$, the right-hand side of (\ref{convergenceeq})
vanishes asymptotically.
{However, if $\llvert \mathbf{C}(h)\rrvert _1=0$ for $\llvert
h\rrvert >h_0$ for some $h_0\in\N
$, setting $l=h_0$ fixed suffices. In this case, the expression on the
right-hand side of (\ref{convergenceeq}) is of faster order $O(n^{-1/2})$.}

As already pointed out by \citet{McMurryPolitis2010}, the tapered
estimator $\widehat{\bolds \Gamma}_{\kappa,l}$ is not guaranteed to
be positive semi-definite or even to be positive definite for finite
sample sizes. However, $\widehat{\bolds \Gamma}_{\kappa,l}$ is at
least ``asymptotically positive definite'' under assumption (A3) and
due to (\ref{convergenceeq}) if
{$1/l+l/\sqrt{n}=o(1)$ holds}. {In the following, we require a
consistent estimator for $\bolds{\Gamma}_{dn}$ which is positive
definite} for all finite
sample sizes to be able to compute its Cholesky decomposition for the
linear process bootstrap scheme that will be introduced in Section~\ref
{bootstrapscheme} below.

{To obtain an estimator of $\bolds{\Gamma}_{dn}$ related to $\widehat{\bolds\Gamma}_{\kappa,l}$ that is assured to be positive definite
for all sample sizes, we construct a modified estimator $\widehat{\bolds\Gamma}_{\kappa,l}^\varepsilon$ in the following. Let
$\widehat\V=\operatorname{diag}(\widehat{\bolds \Gamma}_{dn})$ be the\vspace*{1pt} diagonal
matrix of sample variances, and define $\widehat{\mathbf  R}_{\kappa
,l}=\widehat\V^{-1/2}\widehat{\bolds \Gamma}_{\kappa,l}\widehat \V^{-1/2}$. Now we consider the spectral factorization $\widehat{\mathbf R}_{\kappa,l}=\S\D\S^T$, where $\S$ is an $(dn\times
dn)$ orthogonal matrix and $\D=\operatorname{diag}(r_1,\ldots,r_{dn})$ is the
diagonal matrix containing the eigenvalues of $\widehat{\mathbf
R}_{\kappa,l}$ such that $r_1\geq r_2\geq\cdots\geq r_{dn}$. It is
worth noting that this factorization always exists due to symmetry of
$\widehat{\mathbf  R}_{\kappa,l}$, but that the eigenvalues can be
positive, zero or even negative. Now, define
\begin{equation}
\widehat{\bolds \Gamma}_{\kappa,l}^\varepsilon=\widehat\V ^{1/2}
\widehat{\mathbf  R}_{\kappa,l}^\varepsilon\widehat\V^{1/2} =
\widehat\V^{1/2}\S\D^\varepsilon\S^T\widehat
\V^{1/2}, \label
{gammahatepsilon}
\end{equation}
where $\D^\varepsilon=\operatorname{diag}(r_1^\varepsilon,\ldots,r_{dn}^\varepsilon)$ and
$r_i^\varepsilon=\max(r_i,\varepsilon n^{-\beta})$. Here, $\beta>1/2$ and
$\varepsilon>0$ are user defined constants that ensure the positive
definiteness of $\widehat{\bolds \Gamma}_{\kappa,l}^\varepsilon$.
Contrary to the univariate case discussed in \citet{McMurryPolitis2010}, we propose to adjust the eigenvalues of the (equivariant)
correlation matrix $\widehat{\mathbf  R}_{\kappa,l}$ instead of
$\widehat{\bolds \Gamma}_{\kappa,l}$, which then comes along
without a scaling factor in the definition of $r_i^\varepsilon$. Further,
note that setting $\varepsilon=0$ leads to a positive semi-definite
estimate if $\widehat{\bolds \Gamma}_{\kappa,l}$ is indefinite,
which does not suffice for computing the Cholesky decomposition, and
also that $\widehat{\bolds \Gamma}_{\kappa,l}^\varepsilon$ generally
loses the banded shape of $\widehat{\bolds \Gamma}_{\kappa,l}$.
Theorem~\ref{operatornormconvergence2} below, which extends Theorem~3
in \citet{McMurryPolitis2010}, shows that the modification of the
eigenvalues does affect the convergence results obtained in Theorem
\ref{operatornormconvergence1} just slightly.}

%
\begin{theorem}\label{operatornormconvergence2}
Under the assumptions of Theorem~\ref{operatornormconvergence1}, it holds
\begin{eqnarray} \label{convergenceeq2}
&& \bigl\llVert \rho\bigl(\widehat{\bolds \Gamma}_{\kappa,l}^\varepsilon-
\bolds {\Gamma}_{dn}\bigr)\bigr\rrVert _2\nonumber
\\
&&\qquad \leq
\frac{8Md^2(\lfloor c_\kappa
l\rfloor
+1)}{\sqrt{n}}+4\sum_{h=0}^{\lfloor c_\kappa l\rfloor}
\frac
{\llvert  h\rrvert }{n}\bigl\llvert \mathbf{C}(h)\bigr\rrvert _1+4\sum
_{h=l+1}^{n-1}\bigl\llvert \mathbf{C}(h)\bigr
\rrvert _1
\\
&&\quad\qquad{} +\varepsilon \max_i C_{ii}(0)n^{-\beta}+O \biggl(
\frac{1}{n^{1/2+\beta}} \biggr).\nonumber
\end{eqnarray}
\end{theorem}

In comparison to the upper bound established in Theorem~\ref
{operatornormconvergence1}, two more terms appear on the right-hand
side of (\ref{convergenceeq2}) which do converge as well to zero as
$n$ tends to infinity. Note that the first three summands that the
right-hand sides of (\ref{convergenceeq}) and (\ref
{convergenceeq2}) have in common, remain the leading terms if $\beta
>\frac{1}{2}$.

{We also need convergence and boundedness in operator norm of
quantities related to $\widehat{\bolds \Gamma}_{\kappa,l}^\varepsilon
$. The required results are summarized in the following corollary.}

%
\begin{corollary}\label{operatornormconvergence3}
Under assumptions \textup{(A1)} with $g=0$, \textup{(A2)} and \textup{(A3)}, we have:
\begin{longlist}[(iii)]
\item[(i)] $\rho(\widehat{\bolds \Gamma}_{\kappa,l}^\varepsilon
-\bolds{\Gamma}_{dn})$ and $\rho((\widehat{\bolds \Gamma
}_{\kappa,l}^\varepsilon)^{-1}-\bolds{\Gamma}_{dn}^{-1})$ are terms
of order $O_P(r_{l,n})$, where
\begin{equation}
r_{l,n}=\frac{l}{\sqrt{n}}+\sum_{h=l+1}^\infty
\bigl\llvert \mathbf {C}(h)\bigr\rrvert _1, \label{rn}
\end{equation}
and {$r_{l,n}=o(1)$ if $1/l+l/\sqrt{n}=o(1)$.}
\item[(ii)] $\rho((\widehat{\bolds \Gamma}_{\kappa,l}^\varepsilon
)^{1/2}-\bolds{\Gamma}_{dn}^{1/2})$ and $\rho((\widehat{\bolds
\Gamma}_{\kappa,l}^\varepsilon)^{-1/2}-\bolds{\Gamma}_{dn}^{-1/2})$
are of order $O_P(\log^2(n)\* r_{l,n})$ and {$\log^2(n)r_{l,n}=o(1)$ if
$1/l+\log^2(n)l/\sqrt{n}=o(1)$, and \textup{(A1)} holds for some $g>0$.}
\item[(iii)] {$\rho(\bolds{\Gamma}_{dn})$, $\rho(\bolds{\Gamma
}_{dn}^{-1})$, $\rho(\bolds{\Gamma}_{dn}^{-1/2})$, $\rho(\bolds
{\Gamma}_{dn}^{1/2})$ are bounded from above and below. $\rho(\widehat{\bolds\Gamma}_{\kappa,l}^\varepsilon)$, $\rho((\widehat{\bolds
\Gamma}_{\kappa,l}^\varepsilon)^{-1})$ and $\rho((\widehat{\bolds
\Gamma}_{\kappa,l}^\varepsilon)^{-1/2})$, $\rho((\widehat{\bolds
\Gamma}_{\kappa,l}^\varepsilon)^{1/2})$ are bounded from above and
below (in probability) if $r_{l,n}=o(1)$ and $\log^2(n)r_{l,n}=o(1)$,
respectively.}
\end{longlist}
\end{corollary}

%
\begin{remark}\label{individualbanding}
In Section~\ref{taperedestimator}, we propose to use a global banding
parameter $l$ that down-weights the autocovariance matrices for
increasing lag; that is, the entire matrix $\C(h)$ is multiplied with
the same $\kappa_l(h)$ in (\ref{gammahat}). However, it is possible
to use individual banding parameters $l_{pq}$ for each sequence of
entries $\{C_{pq}(h),h\in\Z\}$, $p,q=1,\ldots,d$ 
as proposed in \citet{Politis2011}, compare also the simulation
section. 
\end{remark}

\subsection{Selection of tuning parameters}\label{selectionsec}
To get a tapered estimate $\widehat{\bolds \Gamma}_{\kappa,l}$ of
the covariance matrix $\bolds{\Gamma}_{dn}$, some parameters have to
be chosen by the practitioner. These are the flat-top taper $\kappa$
and the banding parameter $l$, which are both responsible for the
down-weighting of the empirical autocovariances $\widehat{\mathbf
C}(h)$ with increasing lag $h$.

To select a suitable taper $\kappa$ from the class of functions (\ref
{kappa}), we have to select $c_\kappa\geq1$ and the function $g$
which determine the range of the decay of $\kappa$ to zero for $\llvert  x\rrvert >1$
and its form over this range, respectively. For some examples of
flat-top tapers, compare \citeauthor{Politis2003b} (\citeyear{Politis2003b,Politis2011}). However, the selection
of the banding parameter $l$ appears to be more crucial than choosing
the tapering function $\kappa$ among the family of well-behaved
flat-top kernels as discussed in \citet{Politis2011}. This is comparable to
nonparametric kernel estimation where usually the bandwidth plays a
more important role than the shape of the kernel.

We focus on providing an empirical rule for banding parameter selection
that has already been used in \citet{McMurryPolitis2010} for the
univariate LPB. They make use of an approach primarily proposed in
\citet{Politis2003b} to estimate the bandwidth in spectral density
estimation which has been generalized to the multivariate case in
\citet{Politis2011}. In the following, we adopt this technique based on the
correlogram/cross-correlogram [cf. \citeauthor{Politis2011} (\citeyear{Politis2011}, Section~6)] for our
purposes. Let
\begin{equation}
\widehat R_{jk}(h)=\frac{\widehat C_{jk}(h)}{\sqrt{\widehat
C_{jj}(0)\widehat C_{kk}(0)}}, \qquad j,k=1,\ldots,d
\label
{cross-correlation}
\end{equation}
be the sample (cross-)correlation between the two univariate time
series $(X_{j,t},t\in\Z)$ and $(X_{k,t},t\in\Z)$ at lag $h\in\Z$.
Now, define $\widehat q_{jk}$ as the smallest nonnegative integer such that
\[
\bigl\llvert \widehat R_{jk}(\widehat q_{jk}+h)\bigr
\rrvert <M_0\sqrt{\log _{10}(n)/n}
\]
for $h=1,\ldots,K_n$, where $M_0>0$ is a fixed constant, and $K_n$ is
a positive, nondecreasing integer-valued function of $n$ such that
$K_n=o(\log(n))$. Note that the constant $M_0$ and the form of $K_n$
are the practitioner's choice. As a rule of thumb, we refer to
\citeauthor{Politis2003b} (\citeyear{Politis2003b,Politis2011}) who makes the concrete recommendation $M_0\simeq2$ and
$K_n=\max(5,\sqrt{\log_{10}(n)})$. After having computed $\widehat
q_{jk}$ for all $j,k=1,\ldots,d$, we take
\begin{equation}
\widehat l=\max_{j,k=1,\ldots,d}\widehat q_{jk}
\label{cut-off}
\end{equation}
as a data-driven \emph{global} choice of the banding parameter $l$.
{By setting $\widehat l_{jk}=\widehat q_{jk}$, we get data-driven \emph
{individual} banding parameter choices as discussed in Remark~\ref
{individualbanding}. For theoretical justification of this empirical
selection of a global cut-off point as the maximum over individual
choices and assumptions that lead to successful adaptation, we refer to
Theorem 6.1 in \citet{Politis2011}.}

Note also that for positive definite covariance matrix estimation, that
is, for computing $\widehat{\bolds \Gamma}_{\kappa,l}^\varepsilon$,
one has to select two more parameters $\varepsilon$ and $\beta$, which
have to be nonnegative and might be set equal to one as suggested in
\citet{McMurryPolitis2010}. 

\section{The multivariate linear process bootstrap procedure}\label{bootstrapscheme}

In this section, we describe the multivariate linear process bootstrap
(MLPB) in detail, discuss some modifications and comment on the special
case where the tapered covariance estimator becomes diagonal.

\begin{longlist}
\item[\textit{Step} 1.] Let $\X$ be the $(d\times n)$ data matrix consisting
of $\R^d$-valued time series data $\underline{X}_1,\ldots,\underline
{X}_n$ of sample\vspace*{1pt} size $n$. Compute the centered observations
$\underline{Y}_t=\underline{X}_t-\overline{\underline{X}}$, where
$\overline{\underline{X}}=\frac{1}{n}\sum_{t=1}^n\underline{X}_t$,
let $\Y$ be the corresponding $(d\times n)$ matrix of centered
observations, and define $\underline{Y}=\operatorname{vec}(\Y)$ to be the
$dn$-dimensional vectorized version of $\Y$.

\item[\textit{Step} 2.] Compute $\underline{W}=(\widehat{\bolds \Gamma
}_{\kappa,l}^\varepsilon)^{-1/2}\underline{Y}$, where $(\widehat{\bolds\Gamma}_{\kappa,l}^\varepsilon)^{1/2}$ denotes the lower left
triangular matrix $\L$ of the Cholesky decomposition $\widehat{\bolds\Gamma}_{\kappa,l}^\varepsilon=\L\L^T$.

\item[\textit{Step} 3.] Let $\underline{Z}$ be the standardized version of
$\underline{W}$, that is, $Z_i=\frac{W_i-\overline{W}}{\widehat
\sigma_W}$, $i=1,\ldots,dn$, where $\overline{W}=\frac{1}{dn}\sum_{t=1}^{dn}W_t$ and $\widehat\sigma_W^2=\frac{1}{dn}\sum_{t=1}^{dn}(W_t-\overline{W})^2$.

\item[\textit{Step} 4.] Generate $\underline{Z}^*=(Z_1^*,\ldots,Z_{dn}^*)^T$
by i.i.d. resampling from $\{Z_1,\ldots,\break  Z_{dn}\}$.

\item[\textit{Step} 5.] Compute $\underline{Y}^*=(\widehat{\bolds \Gamma
}_{\kappa,l}^\varepsilon)^{1/2}\underline{Z}^*$, and let $\Y^*$ be the
matrix that is obtained from $\underline{Y}^*$ by putting this vector
column-wise into an $(d\times n)$ matrix, and denote its columns by
$\underline{Y}_1^*,\ldots,\underline{Y}_n^*$.
\end{longlist}

Regarding steps 3 and 4 above and due to the multivariate nature of
the data, it appears to be even more natural to split the
$dn$-dimensional vector $\underline{Z}$ in step 3 above in $n$
sub-vectors, to center and standardize them and to apply i.i.d.
resampling to these vectors to get $\underline{Z}^*$. More precisely,
steps 3 and 4 can be replaced by:

\begin{longlist}
\item[\textit{Step} 3$^\prime$.] Let $\underline{Z}=(\underline
{Z}_1^T,\ldots,\underline{Z}_n^T)^T$ be the standardized version of
$\underline{W}$,
that is, $\underline{Z}_i=\widehat\Sigma_W^{-1/2}(\underline
{W}_i-\overline{\underline{W}})$, where $\underline{W}=(\underline
{W}_1^T,\ldots,\underline{W}_n^T)^T$, $\overline{\underline
{W}}=\frac{1}{n}\sum_{t=1}^{n}\underline{W}_t$ and $\widehat\Sigma
_W=\frac{1}{n}\sum_{t=1}^{n}(\underline{W}_t-\overline{\underline
{W}})(\underline{W}_t-\overline{\underline{W}})^T$.

\item[\textit{Step} 4$^\prime$.] Generate $\underline{Z}^*=(\underline
{Z}_1^*,\ldots,\underline{Z}_n^*)^T$ by i.i.d. resampling from $\{
\underline{Z}_1,\ldots,\break \underline{Z}_{n}\}$.
\end{longlist}

This might preserve more higher order features of the data that are not
captured by $\widehat{\bolds \Gamma}_{\kappa,l}^\varepsilon$.
However, comparative simulations (not reported in the paper) indicate
that the finite sample performance is only slightly affected by this
sub-vector resampling.

%
\begin{remark}\label{remark2}
If\vspace*{-2pt} $0<l<\frac{1}{c_\kappa}$, the banded covariance matrix estimator
$\widehat{\bolds \Gamma}_{\kappa,l}$ (and $\widehat{\bolds \Gamma
}_{\kappa,l}^\varepsilon$ as well) becomes diagonal. In this case and if
steps 3$^\prime$~and~4$^\prime$ are used, the LPB as described
above is equivalent to the classical i.i.d. bootstrap. Here, note the
similarity to the autoregressive sieve bootstrap which boils down to an
i.i.d. bootstrap if the autoregressive order is $p=0$.\vadjust{\goodbreak}
\end{remark}

\section{Bootstrap consistency for fixed time series dimension}\label{asymptoticresults}

\subsection{Sample mean}
In this section, we establish validity of the MLPB for the sample mean.
The following theorem generalizes Theorem 5 of \citet{McMurryPolitis2010} to the multivariate case under somewhat more general conditions.

%
\begin{theorem}\label{validitysamplemean}
{Under assumptions \textup{(A1)} for some $g>0$, \textup{(A2)}, \textup{(A3)}, \textup{(A4)} for $q=4$,
\textup{(A5)} and $1/l+\log^2(n)l/\sqrt{n}=o(1)$}, the MLPB is asymptotically
valid for the sample mean $\overline{\underline{X}}$, that is,
\[
\sup_{\underline{x}\in\R^d}\bigl\llvert P \bigl\{\sqrt{n} (
\overline{\underline {X}}-\underline{\mu} )\leq\underline{x} \bigr\}-P^* \bigl\{
\sqrt{n} \overline{\underline{Y}}^*\leq\underline{x} \bigr\} \bigr\rrvert
=o_P(1)
\]
and $\operatorname{Var}^* (\sqrt{n} \overline{\underline{Y}}^* )=\sum_{h=-\infty}^\infty\C(h)+o_P(1)$, where $\overline{\underline
{Y}}^*=\frac{1}{n}\sum_{t=1}^n \underline{Y}_t^*$.
The short-hand $\underline{x}\leq\underline{y}$ for $\underline
{x},\underline{y}\in\R^d$ is used to denote $x_i\leq y_i$ for all
$i=1,\ldots,d$.
\end{theorem}

\subsection{Kernel spectral density estimates}
Here we prove consistency of the MLPB for kernel spectral density
matrix estimators; this result is novel even in the univariate case.
Let $\mathbf{I}_n(\omega)=\underline{J}_n(\omega)\underline
{J}_n^H(\omega)$ the periodogram matrix, where
\begin{equation}
\underline{J}_n(\omega)=\frac{1}{\sqrt{2\pi n}} \sum
_{t=1}^n \underline{Y}_t
e^{-it\omega} \label{DFT}
\end{equation}
is the discrete Fourier transform (DFT) of $\underline{Y}_1,\ldots
,\underline{Y}_n$, $\underline{Y}_t=\underline{X}_t-\overline
{\underline{X}}$. We define the estimator
\begin{equation}
\widehat{\mathbf{f}}(\omega)=\frac{1}{n} \sum
_{k=-\lfloor\vfrac{n-1}{2}\rfloor}^{\lfloor\sfrac{n}{2}\rfloor} K_b(\omega-
\omega_k) \mathbf{I}_n(\omega_k)
\label{fhat}
\end{equation}
for\vspace*{1pt} the spectral density matrix $\mathbf{f}(\omega)$, where $\lfloor
x\rfloor$ is the integer part of $x\in\R$, $\omega_k= 2\pi\frac
{k}{n}, k=-\lfloor\frac{n-1}{2}\rfloor,\ldots,\lfloor\frac
{n}{2}\rfloor$ are the Fourier frequencies, $b$ is the bandwidth and
$K$ is a symmetric and
{square integrable} kernel\vspace*{1pt} function $K(\cdot)$ that satisfies $\int
K(x)\,dx=2\pi$ and $\int K(u)u^2\,du<\infty$ and we set $K_b(\cdot
)=\frac{1}{b}K(\frac{\cdot}{b})$. Let $\mathbf{I}_n^*(\omega)$ be
the bootstrap analogue of $\mathbf{I}_n(\omega)$ based on $\underline
{Y}_1^*,\ldots,\underline{Y}_n^*$ generated\vspace*{1pt} from the MLPB scheme and
let $\widehat{\mathbf{f}}^*(\omega)$ be the bootstrap analogue of
$\widehat{\mathbf{f}}(\omega)$.

%
\begin{theorem}\label{validitykernelspectral}
{Suppose assumptions \textup{(A1)} with $g\geq0$ specified below, \textup{(A2)}, \textup{(A3)},
\textup{(A4)} for $q=8$ and \textup{(A6)} are satisfied. If $b\rightarrow0$ and
$nb\rightarrow\infty$ such that $nb^5=O(1)$ as well as $1/l+\sqrt
{bl}\log^2(n)+\sqrt{nb}\log^2(n)/l^g$ and $1/k+bk^4+\sqrt{nb}\log
^2(n)/k^g$ for some sequence $k=k(n)$, the MLPB is asymptotically valid
for kernel spectral density estimates $\widehat{\mathbf  f}(\omega)$.
That is, for all $s\in\N$ and arbitrary frequencies $0\leq\omega
_1,\ldots,\omega_s\leq\pi$ (not necessarily Fourier frequencies),
it holds}
\begin{eqnarray*}
&& \sup_{\underline{x}\in\R^{d^2s}}\bigl\llvert P \bigl\{ \bigl(
\sqrt{nb}\bigl(\widehat f_{pq}(\omega_j)-f_{pq}(
\omega _j)\bigr)\dvtx p,q=1,\ldots,d; j=1,\ldots,s \bigr)\leq\underline{x}
\bigr\}
\\
&&\hspace*{28pt}{} -P^* \bigl\{ \bigl(\sqrt{nb}\bigl(\widehat f_{pq}^*(\omega
_j)-\check f_{pq}(\omega_j)\bigr)\dvtx p,q=1,
\ldots,d; j=1,\ldots,s \bigr)\leq\underline{x} \bigr\}\bigr\rrvert
\\
&&\qquad = o_P(1),
\end{eqnarray*}
where $\check f_{pq}(\omega)=\frac{1}{2\pi}\sum_{h=-(n-1)}^{n-1}\kappa_l(h)\widehat C_{pq}(h)e^{-ih\omega}$ and, in
particular,
\begin{eqnarray*}
&& nb\Cov^* \bigl(\widehat f^*_{pq}(\omega),\widehat f^*_{rs}(
\lambda ) \bigr)
\\
&&\qquad =  \bigl(f_{pr}(\omega)\overline{f_{qs}(
\omega)}\delta _{\omega,\lambda}+f_{ps}(\omega)\overline{f_{qr}(
\omega)}\tau _{0,\pi} \bigr)
\frac{1}{2\pi}\int K^2(u)\,du+o_P(1),
\end{eqnarray*}
and $E^*(\widehat f^*_{pq}(\omega))-\check f_{pq}(\omega
)=b^2f_{pq}''(\omega)\frac{1}{4\pi}\int K(u)u^2\,du+p_P(b^2)$, for all
$p,q,r,s=1,\ldots,d$ and all $\omega,\lambda\in[0,\pi]$, respectively.
\end{theorem}

\subsection{Other statistics and LPB-of-blocks bootstrap}
{For statistics $T_n$ contained in the broad class of functions of
generalized means, \citet{JentschPolitis2013} discussed how by using a
preliminary blocking scheme tailor-made for a specific statistic of
interest, the MLPB can be shown to be consistent. This class of
statistics contains estimates $T_n$ of $w(\vartheta)$ with $\vartheta
=E(g(\underline{X}_t,\ldots,\underline{X}_{t+m-1}))$ such that}
\[
T_n=w \Biggl\{\frac{1}{n-m+1}\sum_{t=1}^{n-m+1}g(
\underline {X}_t,\ldots,\underline{X}_{t+m-1}) \Biggr\},
\label{Tn}
\]
for some sufficiently smooth functions $g\dvtx \R^{d\times m}\rightarrow\R
^k$, $w\dvtx \R^k\rightarrow\R$ and fixed \mbox{$m\in\N$}. They propose to
block the data first according to the known function $g$ and to apply
then the (M)LPB to the blocked data. More precisely, the multivariate
LPB-of-blocks bootstrap is as follows:
\begin{longlist}
\item[\textit{Step} 1.] Define $\underline{\widetilde X}_t:=g(\underline
{X}_t,\ldots,\underline{X}_{t+m-1})$, and let $\underline{\widetilde
X}_1,\ldots,\underline{\widetilde X}_{n-m+1}$ be the set of blocked data.

\item[\textit{Step} 2.] Apply\vspace*{1pt} the MLPB scheme of Section~\ref{bootstrapscheme}
to the $k$-dimensional blocked data $\underline{\widetilde X}_1,\ldots
,\underline{\widetilde X}_{n-m+1}$ to get bootstrap observations
$\underline{\widetilde X}_1^*,\ldots,\underline{\widetilde X}_{n-m+1}^*$.

\item[\textit{Step} 3.] Compute $T_n^*=w\{(n-m+1)^{-1}\sum_{t=1}^{n-m+1}\underline{\widetilde X}_t^*\}$.

\item[\textit{Step} 4.] Repeat steps 2 and 3 $B$-times, where $B$ is large, and
approximate the unknown distribution of $\sqrt{n}\{T_n-w(\vartheta)\}
$ by the empirical distribution of $\sqrt{n}\{T_{n,1}^*-T_n\},\ldots,\sqrt{n}\{T_{n,B}^*-T_n\}$.
\end{longlist}

The validity of the multivariate LPB-of-blocks bootstrap for some
statistic $T_n$ can be verified by checking the assumptions of Theorem
\ref{validitysamplemean} for the sample mean of the new process $\{
\underline{\widetilde X}_t,t\in\Z\}$.

\section{Asymptotic results for increasing time series dimension}\label{asymptoticresultsincreasing}

In this section, we consider the case when the time series dimension
$d$ is allowed to increase with the sample size $n$, that is,
$d=d(n)\rightarrow\infty$ as $n\rightarrow\infty$. In particular,
we show consistency of tapered covariance matrix estimates and derive
rates that allow for an asymptotic validity result of the MLPB for the
sample mean in this case.

The recent paper by \citet{CaiRenZhou2012} gives a thorough
discussion of the estimation of Toeplitz covariance matrices for
univariate time series. In their setup, that covers also the
possibility of having multiple datasets from the same data generating
process, \citet{CaiRenZhou2012} establish the optimal rates of convergence
using the two simple flat-top kernels discussed in Section~\ref{taperedestimator}, namely the
truncated (i.e., case of pure banding---no tapering) and the trapezoid
taper. When the strength of dependence is quantified via a
smoothness condition on the spectral density, they show that the
trapezoid is superior to the truncated taper, thus confirming the
intuitive recommendations of \citet{Politis2011}. The asymptotic theory of
\citet{CaiRenZhou2012} allows for increasing number of time series
and increasing sample size, but their framework does not contain the
multivariate time series case, neither for fixed nor for increasing
time series dimension, which will be discussed in this section.

Note that Theorem 1 in \citet{McMurryPolitis2010} for the univariate
case, as well as our Theorem~\ref{operatornormconvergence1} for the
multivariate case of fixed time series dimension, give upper
bounds that are quite sharp, coming within a log-term to the
(Gaussian) optimal rate found in Theorem 2 of \citet{CaiRenZhou2012}.

{Instead of assumptions (A1)--(A5) that have been introduced in
Section~\ref{assumptions} and used in Theorem~\ref
{validitysamplemean} to obtain bootstrap consistency for the sample
mean for fixed dimension $d$, we impose the following conditions on the
sequence of time series process $(\{\underline{X}_t^{(n)},t\in\Z\}
)_{n\in\N}$ of now increasing dimension.}

\subsection{Assumptions}

\begin{longlist}[(A3$^\prime$)]
\item[(A1$^\prime$)] {$(\{\underline{X}_t=(X_{1,t},\ldots
,X_{d(n),t})^T,t\in\Z\})_{n\in\N}$ is a sequence of $\R
^{d(n)}$-valued strictly stationary time series processes with mean
vectors $E(\underline{X}_t)=\underline{\mu}=(\mu_1,\ldots, \mu
_{d(n)})$ and autocovariances $\C(h)=(C_{ij}(h))_{i,j=1,\ldots,d(n)}$
defined as in (\ref{covariance}). Here, $(d(n))_{n\in\N}$ is a
nondecreasing sequence of positive integers such that $d(n)\rightarrow
\infty$ as $n\rightarrow\infty$ and, further, suppose
\[
\sum_{h=-\infty}^\infty \Bigl\{\sup
_{n\in\N} \sup_{i,j=1,\ldots
,d(n)}\llvert h\rrvert
^g\bigl\llvert C_{ij}(h)\bigr\rrvert \Bigr\} <\infty
\]
for some $g\geq0$ to be further specified.
}
\item[(A2$^\prime$)] There exists a constant $M^\prime<\infty$ such
that for all $n\in\N$ and all $h$ with $\llvert  h\rrvert <n$, we have
\[
\sup_{i,j=1,\ldots,d(n)}\Biggl\llVert \sum_{t=1}^{n}(X_{i,t+h}-
\overline {X}_i) (X_{j,t}-\overline{X}_j)-nC_{ij}(h)
\Biggr\rrVert _2\leq M^\prime \sqrt{n}.
\]
\item[(A3$^\prime$)] There exists an $n_0\in\N$ large enough such
that for all $n\geq n_0$ and all $d\geq d_0=d(n_0)$ the eigenvalues
$\lambda_1,\ldots,\lambda_{dn}$ of the $(dn\times dn)$ covariance
matrix $\bolds{\Gamma}_{dn}$ are bounded uniformly away from zero
and from above.
\item[(A4$^\prime$)] Define the sequence of projection operators
$P_k^{(n)}(\underline{X})=E(\underline{X}\llvert \mathcal
{F}_{k}^{(n)})-E(\underline{X}\rrvert \mathcal{F}_{k-1}^{(n)})$ for
$\mathcal{F}_k^{(n)}=\sigma(\underline{X}_t,t\leq k)$, and suppose
\[
\sum_{m=0}^\infty \Bigl\{\sup
_{n\in\N} \sup_{i=1,\ldots,d(n)}\bigl\llVert
P_0^{(n)}X_{i,m}\bigr\rrVert _4 \Bigr
\}<\infty
\]
and
\[
\sup_{n\in\N} \sup
_{i=1,\ldots,d(n)}\llVert \overline{X}_i-\mu_i
\rrVert _4=O\bigl(n^{-1/2}\bigr).
\]
\item[(A5$^\prime$)] {For the sample mean, a Cram\'er--Wold-type CLT
holds true. That is, for any real-valued sequence $\underline
{b}=\underline{b}(d(n))$ of $d(n)$-dimensional vectors with $0<M_1\leq
\llvert \underline{b}(d(n))\rrvert _2^2\leq M_2<\infty$ for all
$n\in\N$ and
$v^2=v_{d(n)}^2=\operatorname{Var}(\sqrt{n} (\underline{b}^T (\overline
{\underline{X}}-\underline{\mu} ) ))$, we have
\[
\sqrt{n} \bigl(\underline{b}^T (\overline {\underline{X}}-
\underline{\mu} ) \bigr)/v\overset{\mathcal {D}} {\longrightarrow}
\mathcal{N}(0,1).
\]
}
\end{longlist}

Assumptions (A1$^\prime$)--(A4$^\prime$) are uniform analogues of
(A1)--(A4), which are required here to tackle the increasing time
series dimension $d$. In particular, (A1$^\prime$) implies
\begin{equation}
\sum_{h=-\infty}^\infty\bigl\llvert \C(h)\bigr
\rrvert _1=O\bigl(d^2\bigr). \label{rate}
\end{equation}
Observe also that the autocovariances $C_{ij}(h)$ are assumed to decay
with increasing lag $h$, that is, in time direction, but they are not
assumed to decay with increasing $\llvert  i-j\rrvert $, that is,
with respect to
increasing time series dimension. Therefore, we have to make use of
square summable sequences in (A5$^\prime$) to get a CLT result. {This
technique has been used, for example, by \citet{LewisReinsel1985} and
\citet{GoncalvesKilian2007} to establish central limit results for the
estimation of an increasing number of autoregressive coefficients.} {A
simple sufficient condition for (A5$^\prime$) is, for example, the
case of $(\{\underline{X}_t=(X_{1,t},\ldots,X_{d(n),t})^T,t\in\Z\}
)_{n\in\N}$ being a sequence of i.i.d. Gaussian processes with
eigenvalues of $E(\underline{X}_t\underline{X}_t^T)$ bounded
uniformly from above and away from zero.}

\subsection{Operator norm convergence for increasing time series dimension}
The following theorem generalizes the results of Theorems~\ref
{operatornormconvergence1} and~\ref{operatornormconvergence2} and of
Corollary~\ref{operatornormconvergence3} to the case where $d=d(n)$ is
allowed to increase with the sample size. In contrast to the case of a
stationary spatial process on the plane $\Z^2$ (where a data matrix is
observed that grows in both directions asymptotically as in our
setting), we do not assume that the autocovariance matrix decays in all
directions. {Therefore, to be able to establish a meaningful theory, we
have to replace (A1)--(A5) by the uniform analogues (A1$^\prime
$)--(A5$^\prime$), and due to (\ref{rate}), an additional factor
$d^2$ turns up in the convergence rate and has to be taken into account.}

%
\begin{theorem}\label{operatornormconvergence3d}
Under assumptions \textup{(A1$^\prime$)} with $g\geq0$ specified below,
\textup{(A2$^\prime$)} and \textup{(A3$^\prime$)}, we have:
%
\begin{longlist}[(iii)]
\item[(i)] $\rho(\widehat{\bolds \Gamma}_{\kappa,l}^\varepsilon
-\bolds{\Gamma}_{dn})$ and $\rho((\widehat{\bolds \Gamma
}_{\kappa,l}^\varepsilon)^{-1}-\bolds{\Gamma}_{dn}^{-1})$ are terms
of order $O_P(d^2\widetilde r_{l,n})$, where
\begin{equation}
\widetilde r_{l,n}=\frac{l}{\sqrt{n}}+\sum
_{h=l+1}^\infty \Bigl\{ \sup_{n\in\N}
\sup_{i,j=1,\ldots,d(n)}\bigl\llvert C_{ij}(h)\bigr\rrvert \Bigr\},
\label{rnd}
\end{equation}
and {$d^2\widetilde r_{l,n}=o(1)$ if $1/l+d^2l/\sqrt{n}+d^2/l^g=o(1)$.}
\item[(ii)] $\rho((\widehat{\bolds \Gamma}_{\kappa,l}^\varepsilon
)^{1/2}-\bolds{\Gamma}_{dn}^{1/2})$ and $\rho((\widehat{\bolds
\Gamma}_{\kappa,l}^\varepsilon)^{-1/2}-\bolds{\Gamma}_{dn}^{-1/2})$
are both terms of order\vspace*{1pt} $O_P(\log^2(dn)d^2\widetilde r_{l,n})$ and
{$\log^2(dn)d^2\widetilde r_{l,n}=o(1)$ if $1/l+\log^2(dn)d^2l/\sqrt
{n}+\log^2(dn)d^2/l^g=o(1)$.}
\item[(iii)] {$\rho(\bolds{\Gamma}_{dn})$,\vspace*{1pt} $\rho(\bolds{\Gamma
}_{dn}^{-1})$, $\rho(\bolds{\Gamma}_{dn}^{-1/2})$ and $\rho
(\bolds{\Gamma}_{dn}^{1/2})$ are bounded from above and below. $\rho
(\widehat{\bolds \Gamma}_{\kappa,l}^\varepsilon)$ and $\rho((\widehat{\bolds\Gamma}_{\kappa,l}^\varepsilon)^{-1})$ as well as $\rho
((\widehat{\bolds \Gamma}_{\kappa,l}^\varepsilon)^{-1/2})$ and $\rho
((\widehat{\bolds \Gamma}_{\kappa,l}^\varepsilon)^{1/2})$ are boun\-ded
from above and below in probability if $d^2\widetilde r_{l,n}=o(1)$ and
$\log^2(dn)d^2\times\break \widetilde r_{l,n}=o(1)$, respectively.}
\end{longlist}
\end{theorem}

{The\vspace*{2pt} required rates for the banding parameter $l$ and the time series
dimension $d$ to get operator norm consistency $\rho(\widehat{\bolds
\Gamma}_{\kappa,l}^\varepsilon-\bolds{\Gamma}_{dn})=o_P(1)$ can be
interpreted nicely. If $g$ is chosen to be large enough, $d^2l/\sqrt
{n}$ becomes the leading term, and there is a trade-off between
capturing more dependence of the time series in time direction (large
$l$) and growing dimension of the time series in cross-sectional
direction (large $d$). }

\subsection{Bootstrap validity for increasing time series dimension}
The subsequent theorem is a Cram\'er--Wold-type generalization of
Theorem~\ref{validitysamplemean} to the case where $d=d(n)$ is allowed
to grow at an appropriate rate with the sample size. {To tackle the
increasing time series dimension and to prove such a CLT result, we
have to make use of appropriate sequences of square summable vectors
$\underline{b}=\underline{b}(d(n))$ as described in (A5$^\prime$) above.}

%
\begin{theorem}\label{validitysamplemeand}
{Under assumptions \textup{(A1$^\prime$)} with $g\geq0$ specified below,
\textup{(A2$^\prime$)}, \textup{(A3$^\prime$)}, \textup{(A4$^\prime$)} for $q=4$, \textup{(A5$^\prime$)}
as well as $1/l+\log^2(dn)d^2l/\sqrt{n}+\break \log^2(dn)d^2/l^g=o(1)$
and $1/k+d^5k^4/n+\log^2(dn)d^2/k^g=o(1)$ for some sequence $k=k(n)$},
the MLPB is asymptotically valid for the sample mean $\overline
{\underline{X}}$. {That is, for any real-valued sequence $\underline
{b}=\underline{b}(d(n))$ of $d(n)$-dimensional vectors with $0<M_1\leq
\llvert \underline{b}(d(n))\rrvert _2^2\leq M_2<\infty$ for all
$n\in\N$ and
$\widehat v^2=\widehat v_{d(n)}^2=\operatorname{Var}^*(\sqrt{n}(\underline
{b}^T\overline{\underline{Y}}^*))$, we have
\[
\sup_{x\in\R}\bigl\llvert P \bigl\{\sqrt{n}
\bigl(\underline{b}^T (\overline{\underline {X}}-
\underline{\mu} ) \bigr)/v\leq x \bigr\}-P^* \bigl\{ \sqrt{n} \bigl(
\underline{b}^T\overline{\underline{Y}}^* \bigr)/\widehat v\leq x
\bigr\}\bigr\rrvert =o_P(1)
\]
and $\llvert  v^2-\widehat v^2\rrvert =o_P(1)$.}
\end{theorem}

{
\subsection{Reduction of computation time}\label{secreduction}
In practice, the computational requirements can become very demanding
for large $d$ and $n$. In this case, we suggest to split the data
vector $\underline{X}$ in few subsamples $\underline{X}^{(1)},\ldots
,\underline{X}^{(S)}$, say,
and to apply the MLPB scheme to each subsample separately. This
{operation} can be justified by the fact that dependence structure is
distorted only few times. Precisely, we suggest the following procedure:
\begin{longlist}
\item[\textit{Step} 1.] For small $S\in\N$, define $n_{\mathrm{sub}}=\lceil n/S\rceil
$ and $N_{\mathrm{sub}}=dn_{\mathrm{sub}}$ such that $SN_{\mathrm{sub}}\geq N$, and let
$\underline{X}^{(i)}=(\underline{X}_{(i-1)n_{\mathrm{sub}}+1}^T,\ldots
,\underline{X}_{in_{\mathrm{sub}}}^T)^T$, $i=1,\ldots,S$, where\vspace*{1pt} $\underline
{X}^{(S)}$ is filled up with zeros if $SN_{\mathrm{sub}}>N$.

\item[\textit{Step} 2.] Apply the MLPB bootstrap scheme as described in
Section~\ref{bootstrapscheme} separately to the subsamples $\underline
{X}^{(1)},\ldots,\underline{X}^{(S)}$ to get $\underline
{Y}^{(1)*},\ldots,\underline{Y}^{(S)*}$.\vspace*{1pt}

\item[\textit{Step} 3.] Put $\underline{X}^{(1)*},\ldots,\underline
{X}^{(S)*}$ end-to-end together, and discard the last $SN_{\mathrm{sub}}-N$
values to get $\underline{Y}^{*}$ and $\mathbf{Y}^*$.
\end{longlist}

Here, computationally demanding operations as eigenvalue decomposition,
Cholesky decomposition and matrix inversion have to be executed only
for lower-dimensional matrices, such that the algorithm above is
capable to reduce the computation time considerably. Further,\vspace*{1pt} to regain
efficiency, we propose to use the pooled sample mean $\overline
{\underline{X}}$ for centering and $\widehat{\bolds \Gamma}_{\kappa
,l,N_{\mathrm{sub}}}^\varepsilon$ for whitening and re-introducing correlation
structure for \emph{all} subsamples in step 2. Here, $\widehat{\bolds\Gamma}_{\kappa,l,N_{\mathrm{sub}}}^\varepsilon$ is obtained
analogously to (\ref{gammahatepsilon}), but based on the upper-left
$(N_{\mathrm{sub}}\times N_{\mathrm{sub}})$ sub-matrix of $\widehat{\bolds \Gamma
}_{\kappa,l}$.
}

\section{Simulations}\label{secsim}
{
In this section we compare systematically the performance of the
multivariate linear process bootstrap (MLPB) to that of the
vector-autoregressive sieve bootstrap (AR-sieve), the moving block
bootstrap (MBB) and the tapered block bootstrap (TBB) by means of
simulation. In order to make such a comparison, we have chosen a
statistic for which all methods lead to asymptotically correct
approximations. Being interested in the distribution of the sample
mean, we compare the aforementioned bootstrap methods by plotting:
\begin{longlist}[(a)]
\item[(a)] root mean squared errors (RMSE) for estimating the
variances of $\sqrt{n}(\overline{\underline{X}}-\underline{\mu})$ and
\item[(b)] coverage rates (CR) of 95\% bootstrap confidence intervals
for the components of $\underline{\mu}$
\end{longlist}
for two data generating processes (DGPs) and three sample sizes in two
different setups. First, in Section~\ref{secsimtuning}, we compare
the performance of all aforementioned bootstraps with respect to (w.r.t.)
tuning parameter choice. These are the banding parameter $l$ (MLPB),
the autoregressive order $p$ (AR-sieve) and the block length~$s$ (MBB,
TBB). Furthermore, we report RMSE and CR for data-adaptively chosen
tuning parameters to investigate how accurate automatic selection
procedures can work in practice. Second, in Section~\ref
{secsimdimension}, we investigate the effect of the time series
dimension $d$ on the performance of the different bootstrap approaches.

For each case, we have generated $T=500$ time series and $B=500$
bootstrap replications have been used in each step. For (a), the exact
covariance matrix of $\sqrt{n}(\overline{\underline{X}}-\underline
{\mu})$ is estimated by 20,000 Monte Carlo replications. Further, we
use the trapezoidal kernel defined in (\ref{trapezoid}) to taper the
sample covariance matrix for the MLPB and the blocks for the TBB. To
correct the covariance matrix estimator $\widehat{\bolds \Gamma
}_{\kappa,l}$ to be positive definite, if necessary, we set $\varepsilon
=1$ and $\beta=1$ to get $\widehat{\bolds \Gamma}_{\kappa
,l}^\varepsilon$. This choice has already been used by \citet{McMurryPolitis2010} and simulation results (not reported in this paper)
indicate that the performance of the MLPB reacts only slightly to this
choice. We have used the sub-vector resampling scheme, that is, steps
3$'$ and 4$'$ described in Section~\ref{bootstrapscheme}.

Some additional simulation results and a real data application of the
MLPB to the weighted mean of an increasing number of German stock
prices taken from the DAX index can be found in the supplementary
material to this paper [\citet{JentschPolitis2014}]. The R code is
available at \url{http://www.math.ucsd.edu/\textasciitilde
politis/SOFT/function\_MLPB.R}.

\subsection{Bootstrap performance: The effect of tuning parameter choice} \label{secsimtuning} We consider realizations $\underline
{X}_1,\ldots,\underline{X}_n$ of
length $n=100, 200, 500$ from two bivariate ($d=2$) DGPs. Precisely, we study
a first-order vector moving average process
\[
\mbox{VMA(1) model} \qquad\underline{X}_t=\mathbf{A}\underline
{e}_{t-1}+\underline{e}_t
\]
and a first-order vector autoregressive process
\[
\mbox{VAR(1) model} \qquad\underline{X}_t=\mathbf{A}\underline
{X}_{t-1}+\underline{e}_t,
\]
where $\underline{e}_t\sim\mathcal{N}(0,\bolds{\Sigma})$ is a normally distributed i.i.d. white noise process and
\[
\bolds{\Sigma}= \pmatrix{1 & 0.5
\cr
0.5 & 1}\quad\mbox{and}\quad\mathbf{A}=
\pmatrix{0.9 & -0.4
\cr
0 & 0.5}
\]
have been used in all cases. It is worth noting that (asymptotically)
all bootstrap procedures under consideration yield valid approximations
for both models above.
For the VMA(1) model, MLPB is valid for all (sufficiently small)
choices of banding parameters $l\geq1$, but AR-sieve is valid only
asymptotically for $p=p(n)$ tending to infinity at an appropriate rate
with increasing sample size $n$. This relationship of MLPB and AR-sieve
is reversed for the VAR(1) model. For the MBB and the TBB, the block
length has to increase with the sample size for both DGPs.

%
\begin{figure}

\includegraphics{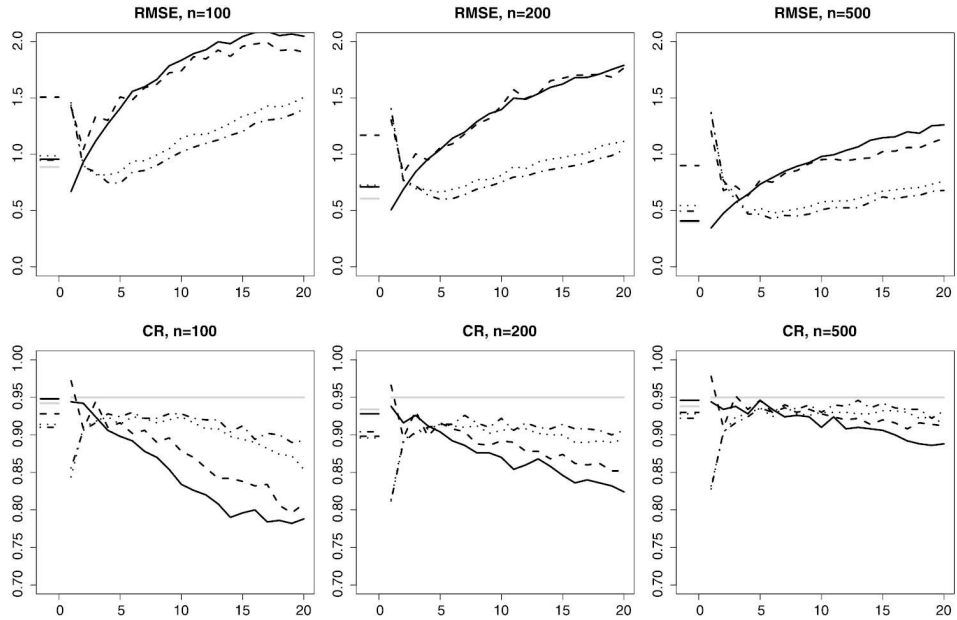}

\caption{{RMSE for estimating $\operatorname{Var}(\sqrt{n}(\overline{X}_1-\mu_1))$
and CR of bootstrap confidence intervals for $\mu_1$ by MLPB (solid),
AR-sieve (dashed), MBB (dotted) and TBB (dash-dotted) are reported vs.
the respective tuning parameters $l,p,s\in\{1,\ldots,20\}$ for the
VMA(1) model with sample size $n\in\{100,200,500\}$. Line segments
indicate results for data-adaptively chosen tuning parameters. MLPB
with individual (grey) and global (black) banding parameter choice are
reported.}} \label{fig2}
\end{figure}

%
\begin{figure}

\includegraphics{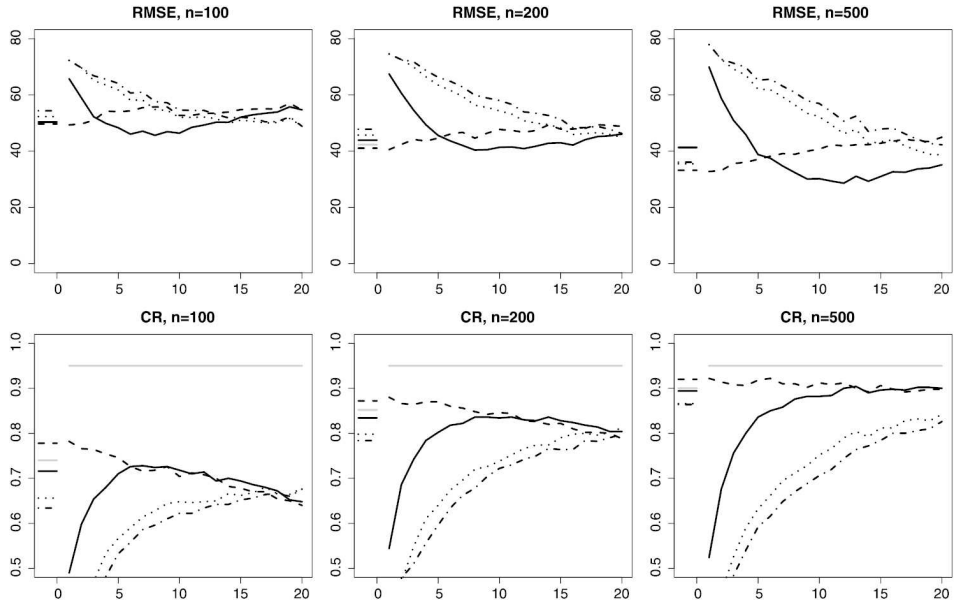}

\caption{As in Figure~\protect\ref{fig2}, but with VAR(1) model.}
\label{fig3}
\end{figure}

In addition to the results for tuning parameters $l,p,s\in\{1,\ldots
,20\}$, we show also RMSE and CR for tuning parameters chosen by
automatic selection procedures in Figures~\ref{fig2} and~\ref{fig3}.
For the MLPB, we report results for data-adaptively chosen global and
individual banding parameters as discussed in Section~\ref
{selectionsec}. For the AR-sieve, the order of the VAR model fitted to
the data has been chosen by using the \textbf{R} routine VAR$(\cdot)$ contained
in the package \textbf{vars} with \textit{lag}.\textit{max}${}=\textit{sqrt}(n/log(n))$. The block
length is chosen by using the \textbf{R} routine \textit{b}.\textit{star}$(\cdot)$ contained in
the package \textbf{np}. In Figures~\ref{fig2} and~\ref{fig3}, we report
only the results corresponding to the first component of the sample
mean, as those for the second component lead qualitatively to the same
results. We show them in the supplementary material, which contains
also corresponding simulation results for a normal white noise DGP.

For data generated by the VMA(1) model, Figure~\ref{fig2} shows that
the MLPB outperforms AR-sieve, MBB and TBB for adequate tuning
parameter choice, that is, $l\approx1$. In this case, the MLPB
generally behaves superiorly, with respect to RMSE and CR, to the other
bootstrap methods for all tuning parameter choices of $p$ and~$s$. This
was not unexpected since, by design, the MLPB can approximate very
efficiently the covariance structure of moving average processes.
Nevertheless, due to the fact that all proposed bootstrap schemes are
valid at least asymptotically, AR-sieve gets rid of its bias with
increasing order $p$, but at the expense of increasing variability and
consequently also increasing RMSE. MLPB with data-adaptively chosen
banding parameter performs quite well, where the individual choice
tends to perform superiorly to the global choice in most cases. In
comparison, MBB and TBB seem to perform quite well for adequate block
length, but they lose in terms of RMSE as well as CR performance if the
block length is chosen automatically.

The data from the VAR(1) model is highly persistent due to the
coefficient $A_{11}=0.9$ near to unity. This leads to autocovariances
that are rather slowly decreasing with increasing lag and,
consequently, to large variances of $\sqrt{n}(\overline{\underline
{X}}-\underline{\mu})$. Figure~\ref{fig3} shows that AR-sieve
outperforms MLPB, MBB and TBB with respect to CR for small AR orders
$p\approx1$. This is to be expected since the underlying VAR(1) model
is captured well by AR-sieve even with finite sample size. But the
picture appears to be different with respect to RMSE. Here, MLPB may
perform superiorly for adequate tuning parameter choice, but this
effect can be explained by the very small variance that compensates its
large bias, in comparison to the AR-sieve (bias and variance not reported here) leading
to a smaller RMSE. This phenomenon is also illustrated by the poor
performance of MLPB with respect to CR for small choices of $l$.
However, more surprising is the rather good performance of the MLPB if
the banding parameter is chosen data-adaptively, where the MLPB appears
to be comparable to the AR-sieve in terms of RMSE and is at least close
with respect to CR. Further, as observed already for the VMA(1) model
in Figure~\ref{fig2}, the individual banding parameter choice
generally tends to outperform the global choice here again. Similarly,
it can be seen here that the performance of AR-sieve worsens with
increasing $p$ at the expense of increasing variability. The block
bootstraps MBB and TBB appear to be clearly inferior to MLPB and
AR-sieve, particularly with respect to CR, but also with respect to
RMSE if tuning parameters are chosen automatically.

\subsection{Bootstrap performance: The effect of larger time series dimension} \label{secsimdimension} We consider $d$-dimensional
realizations $\underline{X}_1,\ldots
,\underline{X}_n$ with $n=100, 200, 500$ from two DGPs of several
dimensions. Precisely, we study first-order vector moving average processes
\[
\mbox{VMA$_d$(1) model} \qquad\underline{X}_t=
\mathbf{A}\underline {e}_{t-1}+\underline{e}_t
\]
and first-order vector autoregressive processes
\[
\mbox{VAR$_d$(1) model} \qquad\underline{X}_t=
\mathbf{A}\underline {X}_{t-1}+\underline{e}_t
\]
of dimension $d\in\{2,\ldots,10\}$, where $\underline{e}_t\sim
\mathcal{N}(0,\bolds{\Sigma}_d)$ is a $d$-dimensional
normally distributed i.i.d. white noise process, and $\bolds{\Sigma
}=(\Sigma_{ij})$ and $A=(A_{ij})$ are such that
\begin{eqnarray*}
\bolds{\Sigma}_{ij}=\cases{ 1, &\quad$i=j$,
\cr
0.5, &\quad$\llvert
i-j\rrvert =1$,
\cr
0, &\quad otherwise} \quad\mbox{and}\quad
\mathbf{A}_{ij}=\cases{ 0.9, &\quad$i=j$, $(i+1)/2\in\N$,
\cr
0.5, &
\quad$i=j$, $i/2\in\N$,
\cr
-0.4, &\quad$i+1=j$,
\cr
0, &\quad otherwise.}
\end{eqnarray*}
Observe that the VMA(1) and VAR(1) models considered in Section~\ref
{secsimtuning} are included in this setup for $d=2$.

In Figures~\ref{fig4} and~\ref{fig5}, we compare the performance of
MLPB, AR-sieve, MBB and TBB for the DGPs above using RMSE and CR
averaged over all $d$ time series coordinates. Precisely, we compute
RMSE individually for the estimates of $\operatorname{Var}(\sqrt{n}(\overline
{X}_i-\mu_i))$, $i=1,\ldots,d$ and plot the averages in the upper
half of Figures~\ref{fig4} and~\ref{fig5}. Similarly, we plot
averages of individually calculated CR of bootstrap confidence
intervals for $\mu_i$, $i=1,\ldots,d$ in the lower halfs. All tuning
parameters are chosen in a data-based and optimal way, as described in
Section~\ref{secsimtuning}, and to reduce computation time, the less
demanding algorithm, as described in Section~\ref{secreduction} with
$S=dn/500$, is used.

%
\begin{figure}

\includegraphics{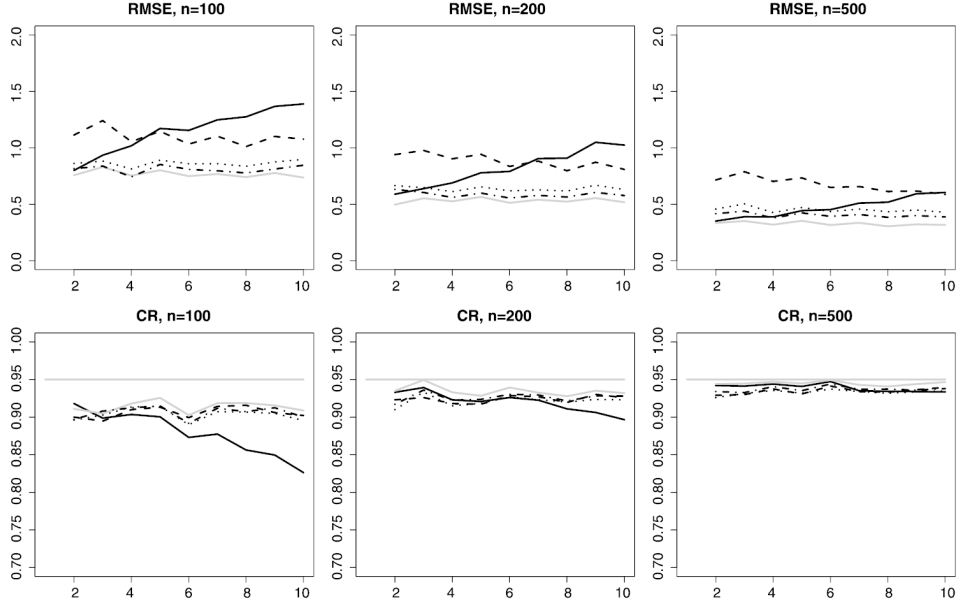}

\caption{Average RMSE for estimating $\operatorname{Var}(\sqrt{n}(\overline
{X}_i-\mu_i))$, $i=1,\ldots,d$ and average CR of bootstrap confidence
intervals for $\mu_i$, $i=1,\ldots,d$, by MLPB (solid), AR-sieve
(dashed), MBB (dotted) and TBB (dash-dotted) with data-based optimal
tuning parameter choices are reported vs. the dimension $d\in\{
2,\ldots,10\}$ for the VMA$_d$(1) model with sample size $n\in\{
100,200,500\}$. MLPB with individual
(grey) and global (black) banding parameter choice are reported.} \label{fig4}
\end{figure}

%
\begin{figure}

\includegraphics{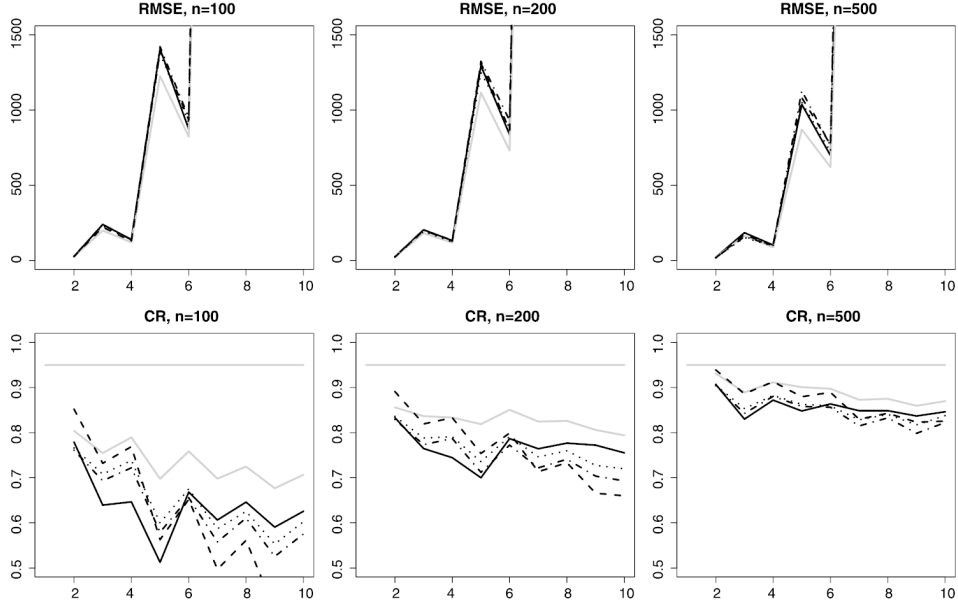}

\caption{As in Figure~\protect\ref{fig4}, but with VAR$_d$(1)
model.}\label{fig5}
\end{figure}

For the VMA(1) DGPs in Figure~\ref{fig4}, the MLPB with individual
banding parameter choice outperforms the other approaches essentially
for all time series dimension under consideration with respect to
averaged RMSE and CR. In particular, larger time series dimensions do
not seem to have a large effect on the performance of all bootstraps
for the VMA(1) DGPs, with the only exception being the MLPB with global
banding parameter choice. In particular, the latter is clearly inferior
in comparison to the MLPB with individually chosen banding parameter,
which might be explained by sparsity of the covariance matrix $\bolds{\Gamma}6
_{dn}$.

In Figure~\ref{fig5}, for the VAR(1) DGPs, the picture is different
from the VMA(1) case above. The influence of larger time series
dimension on RMSE (and less pronounced for CR) performance is much more
pronounced and clearly visible. In particular, the RMSE blows up with
increasing dimension $d$ for all four bootstrap methods, which is due
to the also increasing variance of the process. Note that the zig-zag
shape of the RMSE curves is due to the back and forth switching from
$0.9$ to $0.5$ on the diagonal of $A$. As already observed for the
VMA(1) DGPs, the MLPB with individual banding parameter choice again
performs best over essentially all time series dimensions with respect
to average RMSE and average CR. In particular, MLPB with individual
choice is superior to the global choice. Here, the good performance of
the MLPB is somewhat surprising as the VAR(1) DGPs have rather slowly
decreasing autocovariance structure, where we expected an AR-sieve to
be more suitable.}

\section*{Acknowledgments}
The authors thank Timothy McMurry
for his helpful advice on the univariate case and three anonymous
referees and the Editor who helped to significantly improve the
presentation of the paper.

\begin{supplement}[id=suppA]
\stitle{Additional proofs, simulations and a real data example\\}
\slink[doi]{10.1214/14-AOS1301SUPP} 
\sdatatype{.pdf}
\sfilename{aos1301\_supp.pdf}
\sdescription{In the supplementary material we provide proofs,
additional supporting
simulations and an application of the MLPB to German stock index data.
The supplementary material to this paper is also available online
at \url{http://www.math.ucsd.edu/\textasciitilde politis/PAPER/MLPBsupplement.pdf}.}
\end{supplement}


%

\printaddresses

\begin{thebibliography}{36}
\bibitem[\protect\citeauthoryear{Brillinger}{1981}]{Brillinger1981}
%
\begin{bbook}[mr]
\bauthor{\bsnm{Brillinger},~\bfnm{David~R.}\binits{D.~R.}}
(\byear{1981}).
\btitle{Time Series: Data Analysis and Theory},
\bedition{2nd} ed.
\bpublisher{Holden-Day},
\blocation{Oakland, CA}.
\bid{mr={0595684}}
\end{bbook}
%

\bptok{imsref}%
\endbibitem

\bibitem[\protect\citeauthoryear{Brockwell and Davis}{1988}]{BrockwellDavis1988}
%
\begin{barticle}[mr]
\bauthor{\bsnm{Brockwell},~\bfnm{P.~J.}\binits{P.~J.}} \AND
\bauthor{\bsnm{Davis},~\bfnm{R.~A.}\binits{R.~A.}}
(\byear{1988}).
\btitle{Simple consistent estimation of the coefficients of a linear filter}.
\bjournal{Stochastic Process. Appl.}
\bvolume{28}
\bpages{47--59}.
\bid{doi={10.1016/0304-4149(88)90063-4}, issn={0304-4149}, mr={0936372}}
\end{barticle}
%

\bptok{imsref}%
\endbibitem

\bibitem[\protect\citeauthoryear{Brockwell and
Davis}{1991}]{BrockwellDavis1991}
%
\begin{bbook}[mr]
\bauthor{\bsnm{Brockwell},~\bfnm{Peter~J.}\binits{P.~J.}} \AND
\bauthor{\bsnm{Davis},~\bfnm{Richard~A.}\binits{R.~A.}}
(\byear{1991}).
\btitle{Time Series: Theory and Methods},
\bedition{2nd} ed.
\bseries{Springer Series in Statistics}.
\bpublisher{Springer},
\blocation{New York}.
\bid{doi={10.1007/978-1-4419-0320-4}, mr={1093459}}
\end{bbook}
%

\bptok{imsref}%
\endbibitem

\bibitem[\protect\citeauthoryear{B{\"u}hlmann}{1997}]{Buehlmann1997}
%
\begin{barticle}[mr]
\bauthor{\bsnm{B{\"u}hlmann},~\bfnm{Peter}\binits{P.}}
(\byear{1997}).
\btitle{Sieve bootstrap for time series}.
\bjournal{Bernoulli}
\bvolume{3}
\bpages{123--148}.
\bid{doi={10.2307/3318584}, issn={1350-7265}, mr={1466304}}
\end{barticle}
%

\bptok{imsref}%
\endbibitem

\bibitem[\protect\citeauthoryear{B{\"u}hlmann}{2002}]{Buehlmann}
%
\begin{barticle}[mr]
\bauthor{\bsnm{B{\"u}hlmann},~\bfnm{Peter}\binits{P.}}
(\byear{2002}).
\btitle{Bootstraps for time series}.
\bjournal{Statist. Sci.}
\bvolume{17}
\bpages{52--72}.
\bid{doi={10.1214/ss/1023798998}, issn={0883-4237}, mr={1910074}}
\end{barticle}
%

\bptok{imsref}%
\endbibitem

\bibitem[\protect\citeauthoryear{Cai, Ren and Zhou}{2013}]{CaiRenZhou2012}
%
\begin{barticle}[mr]
\bauthor{\bsnm{Cai},~\bfnm{T.~Tony}\binits{T.~T.}},
\bauthor{\bsnm{Ren},~\bfnm{Zhao}\binits{Z.}} \AND
\bauthor{\bsnm{Zhou},~\bfnm{Harrison~H.}\binits{H.~H.}}
(\byear{2013}).
\btitle{Optimal rates of convergence for estimating {T}oeplitz
covariance matrices}.
\bjournal{Probab. Theory Related Fields}
\bvolume{156}
\bpages{101--143}.
\bid{doi={10.1007/s00440-012-0422-7}, issn={0178-8051}, mr={3055254}}
\end{barticle}
%

\bptok{imsref}%
\endbibitem

\bibitem[\protect\citeauthoryear{Davidson}{1994}]{Davidson1994}
%
\begin{bbook}[mr]
\bauthor{\bsnm{Davidson},~\bfnm{James}\binits{J.}}
(\byear{1994}).
\btitle{Stochastic Limit Theory: An Introduction for Econometricians}.
\bpublisher{Oxford Univ. Press},
\blocation{New York}.
\bid{doi={10.1093/0198774036.001.0001}, mr={1430804}}
\end{bbook}
%

\bptok{imsref}%
\endbibitem

\bibitem[\protect\citeauthoryear{Dedecker et~al.}{2007}]{Dedeckeretal2007}
%
\begin{bbook}[mr]
\bauthor{\bsnm{Dedecker},~\bfnm{J{\'e}r{\^o}me}\binits{J.}},
\bauthor{\bsnm{Doukhan},~\bfnm{Paul}\binits{P.}},
\bauthor{\bsnm{Lang},~\bfnm{Gabriel}\binits{G.}},
\bauthor{\bsnm{Le{\'o}n R.},~\bfnm{Jos{\'e}~Rafael}\binits{J.~R.}},
\bauthor{\bsnm{Louhichi},~\bfnm{Sana}\binits{S.}} \AND
\bauthor{\bsnm{Prieur},~\bfnm{Cl{\'e}mentine}\binits{C.}}
(\byear{2007}).
\btitle{Weak Dependence: With Examples and Applications}.
\bseries{Lecture Notes in Statistics}
\bvolume{190}.
\bpublisher{Springer},
\blocation{New York}.
\bid{mr={2338725}}
\end{bbook}
%

\bptok{imsref}%
\endbibitem

\bibitem[\protect\citeauthoryear{Doukhan}{1994}]{Doukhan1994}
%
\begin{bbook}[mr]
\bauthor{\bsnm{Doukhan},~\bfnm{Paul}\binits{P.}}
(\byear{1994}).
\btitle{Mixing: Properties and Examples}.
\bseries{Lecture Notes in Statistics}
\bvolume{85}.
\bpublisher{Springer},
\blocation{New York}.
\bid{doi={10.1007/978-1-4612-2642-0}, mr={1312160}}
\end{bbook}
%

\bptok{imsref}%
\endbibitem

\bibitem[\protect\citeauthoryear{Gon{\c{c}}alves and
Kilian}{2007}]{GoncalvesKilian2007}
%
\begin{barticle}[mr]
\bauthor{\bsnm{Gon{\c{c}}alves},~\bfnm{S{\'{\i}}lvia}\binits{S.}}
\AND
\bauthor{\bsnm{Kilian},~\bfnm{Lutz}\binits{L.}}
(\byear{2007}).
\btitle{Asymptotic and bootstrap inference for {${\rm AR}(\infty)$}
processes with conditional heteroskedasticity}.
\bjournal{Econometric Rev.}
\bvolume{26}
\bpages{609--641}.
\bid{doi={10.1080/07474930701624462}, issn={0747-4938}, mr={2415671}}
\end{barticle}
%

\bptok{imsref}%
\endbibitem

\bibitem[\protect\citeauthoryear{Hannan}{1970}]{Hannan1970}
%
\begin{bbook}[mr]
\bauthor{\bsnm{Hannan},~\bfnm{E.~J.}\binits{E.~J.}}
(\byear{1970}).
\btitle{Multiple Time Series}.
\bpublisher{Wiley},
\blocation{New York}.
\bid{mr={0279952}}
\end{bbook}
%

\bptok{imsref}%
\endbibitem

\bibitem[\protect\citeauthoryear{H{\"{a}}rdle, Horowitz and
Kreiss}{2003}]{HaerdleHorowitzKreiss2003}
%
\begin{barticle}[auto:parserefs-M02]
\bauthor{\bsnm{H{\"{a}}rdle},~\bfnm{W.}\binits{W.}},
\bauthor{\bsnm{Horowitz},~\bfnm{J.}\binits{J.}} \AND
\bauthor{\bsnm{Kreiss},~\bfnm{J.-P.}\binits{J.-P.}}
(\byear{2003}).
\btitle{Bootstrap methods for time series}.
\bjournal{Int. Stat. Rev.}
\bvolume{71}
\bpages{435--459}.
\end{barticle}
%

\bptok{imsref}%
\endbibitem

\bibitem[\protect\citeauthoryear{Horn and Johnson}{1990}]{HornJohnson1990}
%
\begin{bbook}[mr]
\bauthor{\bsnm{Horn},~\bfnm{Roger~A.}\binits{R.~A.}} \AND
\bauthor{\bsnm{Johnson},~\bfnm{Charles~R.}\binits{C.~R.}}
(\byear{1990}).
\btitle{Matrix Analysis}.
\bpublisher{Cambridge Univ. Press},
\blocation{Cambridge}.
\bnote{Corrected reprint of the 1985 original}.
\bid{mr={1084815}}
\end{bbook}
%

\bptok{imsref}%
\endbibitem

\bibitem[\protect\citeauthoryear{Jentsch and
Politis}{2013}]{JentschPolitis2013}
%
\begin{barticle}[mr]
\bauthor{\bsnm{Jentsch},~\bfnm{Carsten}\binits{C.}} \AND
\bauthor{\bsnm{Politis},~\bfnm{Dimitris~N.}\binits{D.~N.}}
(\byear{2013}).
\btitle{Valid resampling of higher-order statistics using the linear
process bootstrap and autoregressive sieve bootstrap}.
\bjournal{Comm. Statist. Theory Methods}
\bvolume{42}
\bpages{1277--1293}.
\bid{doi={10.1080/03610926.2012.698781}, issn={0361-0926}, mr={3031281}}
\end{barticle}
%

\bptok{imsref}%
\endbibitem

\bibitem[\protect\citeauthoryear{Jentsch and
Politis}{2015}]{JentschPolitis2014}
%
\begin{bmisc}[author]
\bauthor{\bsnm{Jentsch},~\bfnm{C.}\binits{C.}} \AND
\bauthor{\bsnm{Politis},~\bfnm{D.~N.}\binits{D.~N.}}
(\byear{2015}).
\bhowpublished{Supplement to ``Covariance matrix estimation and linear
process bootstrap for multivariate time series of possibly increasing
dimension.''
DOI:\doiurl{10.1214/14-AOS1301SUPP}}.
\bptok{imsref}%
\end{bmisc}
%

\bptok{imsref}%
\endbibitem

\bibitem[\protect\citeauthoryear{Kreiss}{1992}]{Kreiss1992}
%
\begin{bincollection}[mr]
\bauthor{\bsnm{Kreiss},~\bfnm{Jens-Peter}\binits{J.-P.}}
(\byear{1992}).
\btitle{Bootstrap procedures for {${\rm AR}(\infty)$}-processes}.
In \bbooktitle{Bootstrapping and Related Techniques ({T}rier, 1990)}.
\bseries{Lecture Notes in Econom. and Math. Systems}
\bvolume{376}
\bpages{107--113}.
\bpublisher{Springer},
\blocation{Berlin}.
\bid{doi={10.1007/978-3-642-48850-4_14}, mr={1238505}}
\end{bincollection}
%

\bptok{imsref}%
\endbibitem

\bibitem[\protect\citeauthoryear{Kreiss}{1999}]{Kreiss1999}
%
\begin{bmisc}[auto:parserefs-M02]
\bauthor{\bsnm{Kreiss},~\bfnm{J.-P.}\binits{J.-P.}}
(\byear{1999}).
\bhowpublished{Residual and wild bootstrap for infinite order autoregression.
Unpublished manuscript}.
\end{bmisc}
%

\bptok{imsref}%
\endbibitem

\bibitem[\protect\citeauthoryear{Kreiss and
Paparoditis}{2011}]{KreissPaparoditis2011}
%
\begin{barticle}[mr]
\bauthor{\bsnm{Kreiss},~\bfnm{Jens-Peter}\binits{J.-P.}} \AND
\bauthor{\bsnm{Paparoditis},~\bfnm{Efstathios}\binits{E.}}
(\byear{2011}).
\btitle{Bootstrap methods for dependent data: A review}.
\bjournal{J.~Korean Statist. Soc.}
\bvolume{40}
\bpages{357--378}.
\bid{doi={10.1016/j.jkss.2011.08.009}, issn={1226-3192}, mr={2906623}}
\end{barticle}
%

\bptok{imsref}%
\endbibitem

\bibitem[\protect\citeauthoryear{Kreiss, Paparoditis and
Politis}{2011}]{KreissPaparoditisPolitis2011}
%
\begin{barticle}[mr]
\bauthor{\bsnm{Kreiss},~\bfnm{Jens-Peter}\binits{J.-P.}},
\bauthor{\bsnm{Paparoditis},~\bfnm{Efstathios}\binits{E.}} \AND
\bauthor{\bsnm{Politis},~\bfnm{Dimitris~N.}\binits{D.~N.}}
(\byear{2011}).
\btitle{On the range of validity of the autoregressive sieve bootstrap}.
\bjournal{Ann. Statist.}
\bvolume{39}
\bpages{2103--2130}.
\bid{doi={10.1214/11-AOS900}, issn={0090-5364}, mr={2893863}}
\end{barticle}
%

\bptok{imsref}%
\endbibitem

\bibitem[\protect\citeauthoryear{K{\"u}nsch}{1989}]{Kuensch1989}
%
\begin{barticle}[mr]
\bauthor{\bsnm{K{\"u}nsch},~\bfnm{Hans~R.}\binits{H.~R.}}
(\byear{1989}).
\btitle{The jackknife and the bootstrap for general stationary observations}.
\bjournal{Ann. Statist.}
\bvolume{17}
\bpages{1217--1241}.
\bid{doi={10.1214/aos/1176347265}, issn={0090-5364}, mr={1015147}}
\end{barticle}
%

\bptok{imsref}%
\endbibitem

\bibitem[\protect\citeauthoryear{Lahiri}{2003}]{Lahiri2003}
%
\begin{bbook}[mr]
\bauthor{\bsnm{Lahiri},~\bfnm{S.~N.}\binits{S.~N.}}
(\byear{2003}).
\btitle{Resampling Methods for Dependent Data}.
\bpublisher{Springer},
\blocation{New York}.
\bid{doi={10.1007/978-1-4757-3803-2}, mr={2001447}}
\end{bbook}
%

\bptok{imsref}%
\endbibitem

\bibitem[\protect\citeauthoryear{Lewis and
Reinsel}{1985}]{LewisReinsel1985}
%
\begin{barticle}[mr]
\bauthor{\bsnm{Lewis},~\bfnm{Richard}\binits{R.}} \AND
\bauthor{\bsnm{Reinsel},~\bfnm{Gregory~C.}\binits{G.~C.}}
(\byear{1985}).
\btitle{Prediction of multivariate time series by autoregressive model
fitting}.
\bjournal{J. Multivariate Anal.}
\bvolume{16}
\bpages{393--411}.
\bid{doi={10.1016/0047-259X(85)90027-2}, issn={0047-259X}, mr={0793499}}
\end{barticle}
%

\bptok{imsref}%
\endbibitem

\bibitem[\protect\citeauthoryear{Liu and Singh}{1992}]{LiuSingh1992}
%
\begin{bincollection}[mr]
\bauthor{\bsnm{Liu},~\bfnm{Regina~Y.}\binits{R.~Y.}} \AND
\bauthor{\bsnm{Singh},~\bfnm{Kesar}\binits{K.}}
(\byear{1992}).
\btitle{Moving blocks jackknife and bootstrap capture weak dependence}.
In \bbooktitle{Exploring the Limits of Bootstrap ({E}ast {L}ansing,
MI, 1990)}
(\beditor{\bfnm{R.}\binits{R.}~\bsnm{LePage}}
\AND
\beditor{\bfnm{L.}\binits{L.}~\bsnm{Billard}}, eds.).
\bseries{Wiley Ser. Probab. Math. Statist. Probab. Math. Statist.}
\bpages{225--248}.
\bpublisher{Wiley},
\blocation{New York}.
\bid{mr={1197787}}
\end{bincollection}
%

\bptok{imsref}%
\endbibitem

\bibitem[\protect\citeauthoryear{McMurry and Politis}{2010}]{McMurryPolitis2010}
%
\begin{barticle}[mr]
\bauthor{\bsnm{McMurry},~\bfnm{Timothy~L.}\binits{T.~L.}} \AND
\bauthor{\bsnm{Politis},~\bfnm{Dimitris~N.}\binits{D.~N.}}
(\byear{2010}).
\btitle{Banded and tapered estimates for autocovariance matrices and
the linear process bootstrap}.
\bjournal{J. Time Series Anal.}
\bvolume{31}
\bpages{471--482}.
\bnote{Corrigendum: \textit{J.~Time Ser. Anal.} \textbf{33} (2012).}
\bid{doi={10.1111/j.1467-9892.2010.00679.x}, issn={0143-9782}, mr={2732601}}
\bptnote{check related}%
\end{barticle}
%

\bptok{imsref}%
\endbibitem

\bibitem[\protect\citeauthoryear{Mitchell and
Brockwell}{1997}]{BrockwellMitchell1997}
%
\begin{barticle}[mr]
\bauthor{\bsnm{Mitchell},~\bfnm{Heather}\binits{H.}} \AND
\bauthor{\bsnm{Brockwell},~\bfnm{Peter}\binits{P.}}
(\byear{1997}).
\btitle{Estimation of the coefficients of a multivariate linear filter
using the innovations algorithm}.
\bjournal{J. Time Series Anal.}
\bvolume{18}
\bpages{157--179}.
\bid{doi={10.1111/1467-9892.00044}, issn={0143-9782}, mr={1449807}}
\end{barticle}
%

\bptok{imsref}%
\endbibitem

\bibitem[\protect\citeauthoryear{Paparoditis}{2002}]{Paparoditis2002}
%
\begin{bincollection}[mr]
\bauthor{\bsnm{Paparoditis},~\bfnm{Efstathios}\binits{E.}}
(\byear{2002}).
\btitle{Frequency domain bootstrap for time series}.
In \bbooktitle{Empirical Process Techniques for Dependent Data}
\bpages{365--381}.
\bpublisher{Birkh\"auser},
\blocation{Boston, MA}.
\bid{mr={1958790}}
\end{bincollection}
%

\bptok{imsref}%
\endbibitem

\bibitem[\protect\citeauthoryear{Politis}{2001}]{Politis2001}
%
\begin{bincollection}[auto:parserefs-M02]
\bauthor{\bsnm{Politis},~\bfnm{D.~N.}\binits{D.~N.}}
(\byear{2001}).
\btitle{On nonparametric function estimation with infinite-order
flat-top kernels}.
In \bbooktitle{Probability and Statistical Models with Applications}
(\beditor{\bfnm{Ch.}\binits{Ch.}~\bsnm{Charalambides}} et al., eds.)
\bpages{469--483}.
\bpublisher{Chapman \& Hall/CRC},
\blocation{Boca Raton}.
\end{bincollection}
%

\bptok{imsref}%
\endbibitem

\bibitem[\protect\citeauthoryear{Politis}{2003a}]{Politis2003a}
%
\begin{barticle}[mr]
\bauthor{\bsnm{Politis},~\bfnm{Dimitris~N.}\binits{D.~N.}}
(\byear{2003}a).
\btitle{The impact of bootstrap methods on time series analysis:
Silver anniversary of the bootstrap}.
\bjournal{Statist. Sci.}
\bvolume{18}
\bpages{219--230}.
\bid{doi={10.1214/ss/1063994977}, issn={0883-4237}, mr={2026081}}
\end{barticle}
%

\bptok{imsref}%
\endbibitem

\bibitem[\protect\citeauthoryear{Politis}{2003b}]{Politis2003b}
%
\begin{barticle}[mr]
\bauthor{\bsnm{Politis},~\bfnm{Dimitris~N.}\binits{D.~N.}}
(\byear{2003}b).
\btitle{Adaptive bandwidth choice}.
\bjournal{J. Nonparametr. Stat.}
\bvolume{15}
\bpages{517--533}.
\bid{doi={10.1080/10485250310001604659}, issn={1048-5252}, mr={2017485}}
\end{barticle}
%

\bptok{imsref}%
\endbibitem

\bibitem[\protect\citeauthoryear{Politis}{2011}]{Politis2011}
%
\begin{barticle}[mr]
\bauthor{\bsnm{Politis},~\bfnm{Dimitris~N.}\binits{D.~N.}}
(\byear{2011}).
\btitle{Higher-order accurate, positive semidefinite estimation of
large-sample covariance and spectral density matrices}.
\bjournal{Econometric Theory}
\bvolume{27}
\bpages{703--744}.
\bid{doi={10.1017/S0266466610000484}, issn={0266-4666}, mr={2822363}}
\end{barticle}
%

\bptok{imsref}%
\endbibitem

\bibitem[\protect\citeauthoryear{Politis and
Romano}{1992}]{PolitisRomano1992}
%
\begin{barticle}[mr]
\bauthor{\bsnm{Politis},~\bfnm{Dimitris~N.}\binits{D.~N.}} \AND
\bauthor{\bsnm{Romano},~\bfnm{Joseph~P.}\binits{J.~P.}}
(\byear{1992}).
\btitle{A general resampling scheme for triangular arrays of {$\alpha
$}-mixing random variables with application to the problem of spectral
density estimation}.
\bjournal{Ann. Statist.}
\bvolume{20}
\bpages{1985--2007}.
\bid{doi={10.1214/aos/1176348899}, issn={0090-5364}, mr={1193322}}
\end{barticle}
%

\bptok{imsref}%
\endbibitem

\bibitem[\protect\citeauthoryear{Politis and
Romano}{1994}]{PolitisRomano1994}
%
\begin{barticle}[mr]
\bauthor{\bsnm{Politis},~\bfnm{Dimitris~N.}\binits{D.~N.}} \AND
\bauthor{\bsnm{Romano},~\bfnm{Joseph~P.}\binits{J.~P.}}
(\byear{1994}).
\btitle{Limit theorems for weakly dependent {H}ilbert space valued
random variables with application to the stationary bootstrap}.
\bjournal{Statist. Sinica}
\bvolume{4}
\bpages{461--476}.
\bid{issn={1017-0405}, mr={1309424}}
\end{barticle}
%

\bptok{imsref}%
\endbibitem

\bibitem[\protect\citeauthoryear{Politis and
Romano}{1995}]{PolitisRomano1995}
%
\begin{barticle}[mr]
\bauthor{\bsnm{Politis},~\bfnm{Dimitris~N.}\binits{D.~N.}} \AND
\bauthor{\bsnm{Romano},~\bfnm{Joseph~P.}\binits{J.~P.}}
(\byear{1995}).
\btitle{Bias-corrected nonparametric spectral estimation}.
\bjournal{J.~Time Series Anal.}
\bvolume{16}
\bpages{67--103}.
\bid{doi={10.1111/j.1467-9892.1995.tb00223.x}, issn={0143-9782}, mr={1323618}}
\end{barticle}
%

\bptok{imsref}%
\endbibitem

\bibitem[\protect\citeauthoryear{Rissanen and Barbosa}{1969}]{RissanenBarbosa1969}
%
\begin{barticle}[mr]
\bauthor{\bsnm{Rissanen},~\bfnm{J.}\binits{J.}} \AND
\bauthor{\bsnm{Barbosa},~\bfnm{L.}\binits{L.}}
(\byear{1969}).
\btitle{Properties of infinite covariance matrices and stability of
optimum predictors}.
\bjournal{Information Sci.}
\bvolume{1}
\bpages{221--236}.
\bid{mr={0243711}}
\bptnote{check year}%
\end{barticle}
%

\bptok{imsref}%
\endbibitem

\bibitem[\protect\citeauthoryear{Shao}{2010}]{Shao2010}
%
\begin{barticle}[mr]
\bauthor{\bsnm{Shao},~\bfnm{Xiaofeng}\binits{X.}}
(\byear{2010}).
\btitle{The dependent wild bootstrap}.
\bjournal{J. Amer. Statist. Assoc.}
\bvolume{105}
\bpages{218--235}.
\bid{doi={10.1198/jasa.2009.tm08744}, issn={0162-1459}, mr={2656050}}
\end{barticle}
%

\bptok{imsref}%
\endbibitem

\bibitem[\protect\citeauthoryear{Wu and
Pourahmadi}{2009}]{WuPourahmadi2009}
%
\begin{barticle}[mr]
\bauthor{\bsnm{Wu},~\bfnm{Wei~Biao}\binits{W.~B.}} \AND
\bauthor{\bsnm{Pourahmadi},~\bfnm{Mohsen}\binits{M.}}
(\byear{2009}).
\btitle{Banding sample autocovariance matrices of stationary processes}.
\bjournal{Statist. Sinica}
\bvolume{19}
\bpages{1755--1768}.
\bid{issn={1017-0405}, mr={2589209}}
\end{barticle}
%

\bptok{imsref}%
\endbibitem
\end{thebibliography}
\end{document}